# Categorification: tangle invariants and TQFTs

**Catharina Stroppel**

**Abstract**

Based on different views on the Jones polynomial we review representation theoretic categorified link and tangle invariants. We unify them in a common combinatorial framework and connect them via the theory of Soergel bimodules. The influence of these categorifications on the development of 2-representation theory and the interaction between topological invariants and 2-categorical structures is discussed. Finally, we indicate how categorified representations of quantum groups on the one hand and monoidal 2-categories of Soergel bimodules on the other hand might lead to new interesting 4-dimensional TQFTs.

*Dedicated to Igor Frenkel who introduced me to the world of categorification.*

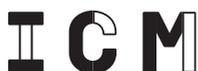

**Introduction**

The study of Topological Quantum Field Theories (TQFTs) is a fruitful interaction between physics and mathematics. The search for interesting TQFTs lead to many developments in mathematical theories which are interesting on their own and also motivates constructions presented in this report. A first mathematical formulation of TQFTs goes back to Atiyah [6], influenced by Segal [130] and Witten [144]. A $d$-(dimensional) TQFT is a symmetric monoidal functor $F$ from a bordism category with objects closed $(d-1)$-dimensional manifolds and morphisms $d$-dimensional bordisms to some symmetric monoidal category, e.g. categories of vector spaces or chain complexes, or more complicated categories.

Representation theory is a good source for TQFTs and monoidal categories. For example, categories of group representations give Dijkgraaf–Witten TQFTs [39] for any $d$. For $d = 3$, representations of quantum groups, i.e. quantised representations of Lie algebras, provide rich and interesting TQFTs due to Reshetikhin–Turaev [119] and Turaev–Viro [141]. They are often viewed as mathematical formulations of Chern–Simons theory [143] or of a form of Ponzano–Regge state sum model from Quantum Gravity [113]. These theories are closely related to (Laurent-)polynomial invariants of knots and links. In Chern–Simons theory [143] for instance the partition function is a 3-manifold invariant, but the expectation values of non-local observables supported on one-dimensional defects, the Wilson lines, give such an invariants of links. The *Jones polynomial* $\mathrm{J}(L)$ arises in this way for the gauge group $SU(2)$. In our setting $\mathrm{J}(L)$ appears as the special $\mathfrak{sl}_2$ example of the RT-*link invariants*. While 3-manifolds are rather well-understood, new 4-TQFTs might help to solve open 4-dimensional (smoothness) problems.

TQFTs provide not only numerical invariants for closed manifolds, but also enjoy good locality properties. In order to compute their values on a complicated closed manifold, one usually *cuts* along lower dimensional submanifolds, assigns data to them and *recombines* this simpler data in a clever way. A cutting principle is common in representation theory: representations are described by decomposing into smaller pieces, by finding simple constituents and their multiplicities in a direct sum decomposition or a Jordan–Hölder filtration and by studying functors on these pieces. Encoding such information combinatorially as character formulas, Poincaré polynomials or Kazhdan–Lusztig polynomials etc. has a long successful history. We call this process *decategorification*.

(Re)Combining or gluing conceptually is a rather new focus motivated partially by TQFTs. In *algebraic categorification* combinatorial data gets interpreted categorically by gluing simple constituents in a predicted way, by realising (Laurent-)polynomials as Poincaré polynomials or Euler characteristics, groups as Grothendieck groups of categories and group homomorphisms as the image of exact functors on a Grothendieck group etc. Moreover, functors are considered in *families* with relations between them described algebraically in terms of (quantised) *Lie algebra* or *Hecke algebra actions*. Classical representation theoretic categories are now viewed as *higher categories* equipped with *categorical actions*. As a byproduct, one obtains new and powerful invariants of links, surfaces or higher dimensional manifolds.

We will summarize and bring together known categorifications of (Laurent-)polynomial link invariants based on the *Jones polynomial* $\mathtt{J}(L)$ and its coloured versions $\mathtt{J}_{\mathrm{col}}(L)$ focusing purely on algebraic-representation theoretic constructions around *Soergel (bi)modules* [132]. The precise meaning of *categorification* will depend on the specific construction:

- Section 2.1: $L$ turns into a complex of graded vector spaces with Euler characteristic $\mathtt{J}(L)$, link cobordisms turn into linear maps;
- Section 2.2: $L$ turns into a complex of bigraded vector spaces whose Euler characteristic is a 2-parameter polynomial which specialises to $\mathtt{J}(L)$;
- Section 2.3: $L$ is viewed as a tangle. Boundaries of tangles turn into graded linear categories, tangles into functors, tangle cobordisms into natural transformations, and $\mathtt{J}(L)$ is the value at a specific element of a map between Grothendieck groups;
- Section 2.4: boundaries of *coloured* tangles turn into graded linear categories of possibly infinite global dimension, tangles into functors, and $\mathtt{J}_{\mathrm{col}}(L)$ is the value at a specific element of a map between *completed* Grothendieck groups.

In Section 1 we set up a framework on the decategorified level with different approaches to the Jones polynomial, all coming from quantum groups and Hecke algebras. It unifies and also stresses the differences of the later categorifications. The material is known, but combined from several sources and carefully adapted. An unusual parameter $\eta$ is introduced in order to fit all the normalisations and categorified theories into one common setup.

In Section 2 the pioneering Khovanov invariant $\mathtt{Kh}$ [80] is described first. Recent advances in categorified link invariants indicate that this theory has interesting topological applications and the chance to provide a 4-TQFT [102]. The second categorification we deal with is the triply graded Khovanov–Rozansky link invariant $\mathtt{KR}$ [82], [77], presented in representation theoretic terms. Its values are Laurent series in $v$ over a polynomial ring in two variables. A three parameter *superpolynomial* invariant of links was predicted on the physics side in [40] and constructed for torus links via refined Chern–Simons theory [3]. Connections to double affine Hecke algebras indicated by the appearance of generalised Verlinde algebras were explored mathematically e.g. in [31], [59]. For torus links, superpolynomials can be matched with the $\mathtt{KR}$ invariants by explicit calculations which substantially use categorified Young projectors. Such projectors were introduced in [34], [126], [53] in the context of categorifications of coloured Reshetikhin–Turaev link invariants [34], [53] and are now important tools in the categorified representation theory of Hecke algebras. The third categorification we describe are the Lie theoretic Mazorchuk–Stroppel–Sussan tangle invariants $\mathtt{MSS}^{\pm}$ [135], [105], [139] which implicitly include the $\mathfrak{sl}_k$ Khovanov–Rozansky link invariant [82] via $\mathtt{MSS}^{-}$. Up to some sign issues which appear when passing from webs to matrix factorisations, the two constructions are even connected by a functor. This follows from the Uniqueness Theorem 2.45, Theorem 2.44 and [98]. The two Lie theoretic constructions $\mathtt{MSS}^{\pm}$ (connected by Koszul duality) go one step further: tensor products of representations of quantum $\mathfrak{gl}_k$ are lifted to categories, the action of quantum $\mathfrak{gl}_k$ to functors and the resulting invariant of tangles has values in the homotopy category of some exact functors. Via a categorified *q-skew Howe duality*, the action of tangles can again be expressed in terms of a quantum group action. This allows to put the construction into the setup from [32], [125]. This is an *axiomatic definition* of categorifications of representations of Lie algebras and provides a conceptual 2-categorical framework, where also uniqueness results are established. Quantum group representations and their tangle invariants are then finally turned into 2-*representations of categorified quantum groups* introduced by Khovanov–Lauda [81] and Rouquier [125].

The MSS$^+$-invariant is equivalent to the quantum $\mathfrak{sl}_k$ version from [142] which deals more generally with quantum invariants for *any* reductive Lie algebra. In comparison, [142] *is defined to fit into* and *substantially further developed* the framework of 2-representations, whereas MSS$^+$ *leads naturally to* and *motivates* the framework of 2-representations. In the MSS-theory one should not expect generalisations to other Lie algebra as in [142], but rather to braid groups of general type and to categorified representation theory of certain quantum symmetric pairs via [41] and probably to Khovanov–Rozansky invariants for orthogonal Lie algebras. Important for us is that MSS$^\pm$ *directly* connects to Soergel bimodules, a possible source for an intriguing connection between KR and categorified coloured tangle invariants described as the fourth example. This connection and the complicated combinatorics of categorified coloured link and tangle invariants [54], [34], [142], [138] needs still to be explored.

In Section 3 we return to our motivation: we indicate two (partially conjectural) new approaches towards potentially rich 4-TQFTs, one via categorified representations of quantum groups, the other via semistrict monoidal 2-categories of Soergel bimodules.

**Conventions:** We denote $\mathbb{N} = \mathbb{Z}_{>0}$, $\mathbb{N}_0 = \mathbb{Z}_{\geq 0}$. We fix $\mathbb{C}$ as ground field. Let $S_n$ be the symmetric group on *n* letters with standard generators $s_i = (i, i+1)$, $1 \leq i \leq n-1$, and length function $\ell$. For a variable $v$ and $a \in \mathbb{Z}$ let $[a] := \frac{v^a - v^{-a}}{v - v^{-1}} = v^{a-1} + v^{a-3} + \cdots + v^{1-a} \in \mathbb{Z}[v^{\pm 1}]$ be the *v-quantum number*, a Laurent polynomial in $v$. By *graded* we mean $\mathbb{Z}$-graded and $\langle i \rangle$ denotes the shift up in the grading, i.e. $(M\langle i \rangle)_n = (\langle i \rangle M)_n = M_{n-i}$. Similarly we write $[i]$ for the shift of complexes by $i$ in direction of the differential. When displaying complexes, we indicate the homological degree zero by putting a box around the component. For an additive category $\mathcal{A}$, we denote by $K^b(\mathcal{A})$ the homotopy category of bounded complexes in $\mathcal{A}$. When describing morphisms or functors diagrammatically we read from bottom to top, and composition is vertical stacking, whereas a monoidal product $\otimes$ is denoted by horizontal juxtaposition, and identities are usually displayed by a vertical strand.

### 1. Four approaches to the Jones polynomial

We summarize four similar, but different, algebraic approaches to knot or link invariants giving rise to the $\mathbb{Z}[v^{\pm 1}]$-valued Jones polynomial. These approaches will later be connected with four theories in the context of categorification. The third and fourth, RT and wRT, are more involved and cover also tangles (a common generalisation of links and braids). The first one is best for computations, but the passage to tangles requires extra adjustments like the use of skein algebras. The second one does not cover tangles at all, but is probably the most intuitive approach for categorifications. It works with link closures instead of planar projections of links. In the following, $v$ denotes a (generic) variable.

**I. Kauffman bracket of links.** We fix an orientation of $\mathbb{R}^3$ and consider oriented knots or links $L$ in $\mathbb{R}^3$. Following Kauffman [73] we first ignore the orientation and assign to any generic, i.e. with no triple intersections, no tangencies and no cusps, planar projection $D$ of $L$, the *Kauffman bracket* $[\![D]\!] \in \mathbb{Z}[v^{\pm 1}]$. It is characterised by the *multiplicativity property* $[\![D_1 \sqcup D_2]\!] := [\![D_1]\!][\![D_2]\!]$, i.e. the bracket of a disjoint union is the product of the brackets, and the following *normalisation* and *local smoothing* relation (which removes crossings):

(1.1) $\quad [\![\bigcirc]\!] = v + v^{-1} = [2]$ and $[\![\times]\!] = [\![\smile\frown]\!] - v[\![\,|\,|\,]\!]$ respectively.

The assignment $D \mapsto \mathtt{J}(D) := (-1)^{n_-(D)} v^{n_+(D) - 2n_-(D)} [\![D]\!] \in \mathbb{Z}[v^{\pm 1}]$, where $n_\pm(D)$ denotes the number of positive respectively negative crossings in $D$, defines then an invariant of oriented links, the *Jones polynomial* $\mathtt{J}(D)$. It fulfils the skein relation, with $\mathbf{a} = v^2$ and $\mathtt{P}(D) = \mathtt{J}(D)$,

(1.2) $$\mathbf{a}\mathtt{P}\left(\nearrow\!\!\!\nwarrow\right) - \mathbf{a}^{-1}\mathtt{P}\left(\nwarrow\!\!\!\nearrow\right) = (v - v^{-1})\mathtt{P}\left(\uparrow\uparrow\right).$$

**Example 1.1.** *For the Hopf link diagram $D = $ ⧖ the Kauffman bracket has the value*

(1.3) $$[\![\,0\,]\!] - v[\![\,\ominus\,]\!] - v[\![\,\ominus\,]\!] + v^2[\![\,\ominus\ominus\,]\!] = \underbrace{[2]^2 - v[2] - v[2] + v^2[2]^2 = v[4]}_{\Rightarrow \text{ Jones polynomial } \mathtt{J}(D) = v^3[4]}.$$

**II. Closures of braids.** Due to Alexander's theorem [4], every oriented link can be realised as the closure of some upwards oriented braid, i.e. of an element in the usual braid group $\mathrm{Br}_n = \langle \beta_1, \ldots, \beta_{n-1} \rangle$ for some $n$ (see the Hopf link above with $n = 2$). A *Markov trace* Tr with values in some target Inv is a function $\mathrm{Tr}: \coprod_{n \geq 1} \mathrm{Br}_n \to \mathrm{Inv}$ satisfying the *trace condition* $\mathrm{Tr}(\alpha\alpha') = \mathrm{Tr}(\alpha'\alpha)$ and $\mathrm{Tr}(\alpha) = \mathrm{Tr}(\alpha\beta_n^{\pm 1})$ for every $\alpha, \alpha' \in \mathrm{Br}_n$, $n \geq 1$. By Markov's theorem (announced in [103], proved in [18]), Tr induces a well-defined map on isomorphism classes of closures of braids, hence defines an invariant of oriented links. This is a conceptual method to pass from braid invariants to (families) of link invariants. There is an important Markov trace, the *Ocneanu trace* (2.6). Its link invariant is the HOMFLY-PT polynomial $P(L)(v, \mathbf{a}) \in \mathbb{C}(v)[\mathbf{a}^{\pm 1}]$ introduced in [55], [114]. It satisfies (1.2) and $P(L)(v, v^2) = \mathtt{J}(v)$, see Remark 2.11.

**III. Quantum Invariants.** The Jones polynomial of oriented links also arises as the *(Witten)–Reshetikhin–Turaev (=RT) invariant* [120] associated with quantum $\mathfrak{gl}_k$ in the special case $k = 2$. An oriented link is a special case of an oriented tangle, i.e. a disjoint embedding of finitely many arcs and circles into $\mathbb{R}^2 \times [0, 1]$ (sending endpoints of arcs to boundary points) modulo ambient isotopy fixing the boundary points. The RT–invariant assigns to each generic horizontal cut of a tangle a tensor product of modules for the quantum group $U_v(\mathfrak{gl}_k)$, and to each tangle a homomorphism in a consistent way, see Overview 1.

To make this more precise we consider a tangle as a morphism in the monoidal category $\mathcal{T}an$ of oriented tangles with stacking as composition and juxtaposition as tensor product. The *quantum group* $U_v(\mathfrak{gl}_k)$ is a deformation of the universal enveloping Hopf algebra of $\mathfrak{gl}_k$ and is often described as the $\mathbb{C}(v)$-algebra with generators $E_i, F_i, D_j^{\pm 1}$, $1 \leq i \leq k - 1$, $1 \leq j \leq k$ "quantising" the usual matrix units $E_{i,i+1}$, $E_{i+1,i}$, $\pm E_{i,i}$, modulo quantised Serre relations, see e.g. [26] for a definition. It is still a Hopf algebra but now with an interesting

**Overview 1** Hopf link: The RT invariant (see 1. III) and its web version (see 1. IV) with categorifications (see 2.3)

| | RT | wRT | MSS$^-$ | MSS$^+$ |
|---|---|---|---|---|
| | | ($\beta = -v^k H_1^{-1}$) | for $W = V_-$ : | for $W = V_+$ : |
| | $\mathbb{C}(v)$ | $\bigwedge^k W \otimes \bigwedge^k W$ | $\hat{\mathcal{O}}_\mathfrak{c}^{(k,k)}(\mathfrak{g}_{2k})$ | $\hat{\mathcal{O}}_{(k,k)}^\mathfrak{c}(\mathfrak{g}_{2k})$ |
| | ↑ ev | ↑ $\wedge_{1,k-1}^k$ | ↑ $\hat{Z}$ | ↑ |
| | $W \otimes W^*$ | $W \otimes \bigwedge^k W \otimes \bigwedge^{k-1} W$ | $\hat{\mathcal{O}}_\mathfrak{c}^{1,k,k-1}(\mathfrak{g}_{2k})$ | $\hat{\mathcal{O}}_{k-1,k,1}^\mathfrak{c}(\mathfrak{g}_{2k})$ |
| | ↑ id ⊗ ev ⊗ id | ↑ id ⊗ $\wedge_{1,k-1}^k$ ⊗ id | ↑ $\hat{Z}$ | ↑ |
| | $W \otimes W \otimes W^* \otimes W^*$ | $W \otimes W \otimes \bigwedge^{k-1} W \otimes \bigwedge^{k-1} W$ | $\hat{\mathcal{O}}_\mathfrak{c}^{1,1,k-1,k-1}(\mathfrak{g}_{2k})$ | $\hat{\mathcal{O}}_{k-1,k-1,1,1}^\mathfrak{c}(\mathfrak{g}_{2k})$ |
| | ↑ $\beta$ | ↑ $\beta$ | ↑ (2.16) | ↑ |
| | $W \otimes W \otimes W^* \otimes W^*$ | $W \otimes W \otimes \bigwedge^{k-1} W \otimes \bigwedge^{k-1} W$ | $\hat{\mathcal{O}}_\mathfrak{c}^{1,1,k-1,k-1}(\mathfrak{g}_{2k})$ | $\hat{\mathcal{O}}_{k-1,k-1,1,1}^\mathfrak{c}(\mathfrak{g}_{2k})$ |
| | ↑ $\beta$ | ↑ $\beta$ | ↑ (2.16) | ↑ |
| | $W \otimes W \otimes W^* \otimes W^*$ | $W \otimes W \otimes \bigwedge^{k-1} W \otimes \bigwedge^{k-1} W$ | $\hat{\mathcal{O}}_\mathfrak{c}^{1,1,k-1,k-1}(\mathfrak{g}_{2k})$ | $\hat{\mathcal{O}}_{k-1,k-1,1,1}^\mathfrak{c}(\mathfrak{g}_{2k})$ |
| | ↑ id ⊗ coev ⊗ id | ↑ id ⊗ $\vee_k^{1,k-1}$ ⊗ id | ↑ $\hat{\text{incl}}$ | ↑ |
| | $W \otimes W^*$ | $W \otimes \bigwedge^k W \otimes \bigwedge^{k-1} W$ | $\hat{\mathcal{O}}_\mathfrak{c}^{1,k,k-1}(\mathfrak{g}_{2k})$ | $\hat{\mathcal{O}}_{k-1,k,1}^\mathfrak{c}(\mathfrak{g}_{2k})$ |
| | ↑ coev | ↑ $\vee_k^{1,k-1}$ | ↑ $\hat{\text{incl}}$ | ↑ |
| | $\mathbb{C}(v)$ | $\bigwedge^k W \otimes \bigwedge^k W$ | $\hat{\mathcal{O}}_\mathfrak{c}^{(k,k)}(\mathfrak{g}_{2k})$ | $\hat{\mathcal{O}}_{(k,k)}^\mathfrak{c}(\mathfrak{g}_{2k})$ |

non-cocommutative comultiplication (due to the appearance of some $D_j$'s):

(1.4) $\Delta(E_i) = E_i \otimes 1 + D_i D_{i+1}^{-1} \otimes E_i, \quad \Delta(F_i) = 1 \otimes F_i + F_i \otimes D_i^{-1} D_{i+1}, \quad \Delta(D_j^{\pm 1}) = D_j^{\pm 1} \otimes D_j^{\pm 1}$.

Every finite-dimensional representation of $\mathfrak{gl}_k$ quantises to a $U_v(\mathfrak{gl}_k)$-module. As often in quantum algebra, there are different choices for such a quantisation, but for irreducible representations they only differ by a one-dimensional twist. We encode the choice by a function $\eta : \{1, \ldots, k\} \to \{\pm 1\}$ such that the spectrum of $D_j$ is contained in $\eta(j)v^\mathbb{Z}$. To capture different normalisations of link invariants, we at least need to consider the additive monoidal subcategory generated by the irreducibles corresponding to constant $\eta = \pm 1$. These signs, although annoying in practice, often have a deeper meaning in categorifications.

**Example 1.2.** *The quantisation $V_\pm = V_{\pm,\mathfrak{gl}_k}$ of the natural representation of $\mathfrak{gl}_k$ for the constant functions $\eta = \pm 1$ can be realised as the $k$-dimensional $\mathbb{C}(v)$-vector space with basis $e_r$, $1 \leq r \leq k$, and the following $U_v(\mathfrak{gl}_k)$-actions*

(1.5) $\quad V_+ : \quad E_i e_r = \delta_{i,r} e_{r+1}, \quad F_i e_{r+1} = \delta_{i,r} e_r, \quad D_j e_r = v^{\delta_{j,r}} e_r,$
$\quad V_- : \quad E_i e_{r+1} = \delta_{i,r} e_r, \quad F_i e_r = -\delta_{i,r} e_{r+1}, \quad D_j e_r = -v^{\delta_{j,r}} e_r.$

A crucial observation behind the invention of quantum groups was that the permutation action of the symmetric group on tensor products of representations quantises (i.e. lifts) to an action of the braid group. A modern formulation is that $\mathcal{R}ep_k$ is (non-symmetric!) *braided monoidal*. In particular, $\text{Br}_n$ acts on $V_\eta^{\otimes n}$ by $U_v(\mathfrak{gl}_k)$-homomorphisms. Explicitly, $\beta_i$ acts on the $i$-th and $(i+1)$-th tensor factor of $V_\eta^{\otimes n}$ for constant $\eta = \pm 1$ as

(1.6) $\quad H_i : \quad e_a \otimes e_b \mapsto \begin{cases} e_b \otimes e_a & \text{if } a > b, \\ e_b \otimes e_a + (v^{-1} - v) e_a \otimes e_b & \text{if } a < b, \\ \gamma e_a \otimes e_a & \text{if } a = b, \end{cases} \quad \text{with } \gamma := \eta v^{-\eta}.$

These actions factor through $\mathbb{C}(v) \otimes_{\mathbb{Z}[v^{\pm 1}]} \mathbb{H}_n$, where $\mathbb{H}_n$ is the *Hecke algebra*. We define $\mathbb{H}_n$ as the $\mathbb{Z}[v^{\pm 1}]$-algebra quotient of the group algebra $\mathbb{Z}[v^{\pm 1}][\mathrm{Br}_n]$ by the following *quadratic relation*, and denote the image of $\beta_i$ in $\mathbb{H}_n$ or $\mathbb{C}(v) \otimes_{\mathbb{Z}[v^{\pm 1}]} \mathbb{H}_n$ by abuse of notation also $H_i$:

$$(1.7) \quad -\beta_i + \beta_i^{-1} = v - v^{-1} \quad \text{or equivalently} \quad (\beta_i + v)(\beta_i - v^{-1}) = 0.$$

Set $W = V_\eta$. Following [120], the *(Witten)–Reshetikhin–Turaev functor* associated with $W$ is now a monoidal functor $\mathrm{RT} = \mathrm{RT}_W : \mathcal{T}an \to \mathcal{R}ep_k$. It sends an oriented tangle $t$ with, say, $m$ endpoints at the bottom and $n$ endpoints at the top to a $U_v(\mathfrak{gl}_k)$-homomorphism

$$(1.8) \quad \mathrm{RT}(t) : \quad W^{\epsilon_1} \otimes \cdots \otimes W^{\epsilon_m} \longrightarrow W^{\epsilon'_1} \otimes \cdots \otimes W^{\epsilon'_n},$$

where $W^{\epsilon_i} = W$ or $W^{\epsilon_i} = W^*$ respectively, depending whether the $i$th-strand on the bottom of $t$ is oriented up- or downwards, similarly for the top using $W^{\epsilon'_i}$. Now $\mathrm{RT}_W$ is determined by the values on the *elementary* tangles 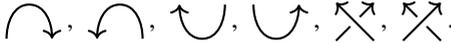. These are sent to the corresponding evaluations, coevaluations and to the morphism $-v^{-k} H_i$ from (1.6) and its inverse $-v^k H_i^{-1}$, respectively. To compute $\mathrm{RT}(t)$ one first reads a chosen generic tangle projection from bottom to top as a vertical composition of *basic tangle diagrams*, i.e. of those which differ from an elementary one just by adding some strands to the left or right, see Overview 1. Each basic tangle diagram is sent to the value of the elementary diagram with identities tensored on the left or right. Finally, $\mathrm{RT}(t)$ is the composition of the values of the basic tangle diagrams. If $t$ is a link, the result is an endomorphism $f$ of $\mathbb{C}(v)$. Evaluating at 1 gives $f(1) \in \mathbb{C}(v)$ which equals $\mathrm{J}(L)$ in case $k = 2$, $W = V_-$. The Hecke relation (1.7) implies the $\mathfrak{gl}_k$-skein relation (1.2) with $\mathbf{a} = v^k$. The RT constructions work for arbitrary reductive Lie algebras, not only $\mathfrak{gl}_k$, and thus provide several families of tangle invariants.

**Remark 1.3.** The (unusual) choice of $V_-$ over $V_+$ has the advantage that the unknot has the value $[k] \in \mathbb{N}[v^{\pm 1}]$ (with nonnegative coefficients!) instead of $(-1)^{k-1}[k]$.

**Remark 1.4.** The construction also works if we pick an irreducible representation for each component of $t$ and gives the *coloured* RT tangle invariant of framed tangles from [120] and the *coloured Jones polynomial* for links if $k = 2$. Colouring only with $V_\eta$ makes life easier, e.g. one can avoid framings and all constructions are defined over $\mathbb{Z}[v^{\pm 1}]$, see Example 2.53.

**IV. Webs and spin networks.** A fourth way to get the Jones polynomial is via webs or spin networks and their evaluations. Following Penrose [112], a web is a certain labelled graph built from trivalent vertices, where a vertex may be interpreted as an event in which either a single unit splits into two or two units collide and join into a single one. More precisely, let $\eta \in \{\pm 1\}$. The *universal* $\mathfrak{gl}$-*web category* is the monoidal $\mathbb{C}(v)$-linear category $\mathcal{W}^\eta$ which is the linear additive closure of the strict monoidal category generated (as monoidal category) by the set of objects $\mathbb{N}$, and on the level of morphisms by diagrams

$$(1.9) \quad \underset{a+b}{\overset{a \quad b}{\bigvee}} \text{ (from } a+b \text{ to } a \otimes b\text{)}, \qquad \underset{a \quad b}{\overset{a+b}{\bigwedge}} \text{ (from } a \otimes b \text{ to } a+b\text{)},$$

modulo the following *associativity* and *coassociativity* relations and *thin square switches*

$$
\begin{array}{c}
\text{(diagrams: associativity, coassociativity, and square switch relations)}
\end{array}
$$

By convention, thin square switches include the *digon removals* (1.10) as degenerate $a = 0$ (or $b = 0$) cases. Together with (co)associativity one obtains *thick square switches* expressing $\boxminus(a,r,b) - \boxminus(a,r,b)$ from (1.11) in terms of a sum over thinner squares, see e.g. [28].

$$
(1.10)\quad a{-}1 \bigcirc 1 = (-\eta)^{a-1}[a] \bigg| = 1 \bigcirc a{-}1 \qquad (1.11)\quad \boxminus(a,r,b) := \cdots \quad \boxminus(a,r,b) := \cdots
$$

An object in $\mathcal{W}^\eta$ is just a finite sequence of nonnegative integers including the empty sequence as tensor unit; a morphism is a linear combination of *webs* obtained by gluing horizontally and vertically the generating pieces (1.9) with identities drawn by vertical lines. For fixed $k \in \mathbb{N}$, let $\mathcal{W}^\eta_k$ be the quotient of $\mathcal{W}^\eta$ by all morphisms factoring through an object involving a number $> k$. We will see that this category provides a concrete graphical presentation of the monoidal category of $U_v(\mathfrak{gl}_k)$-modules generated by quantisations of the fundamental representations of $\mathfrak{gl}_k$. Thus it continues pioneering works on graphical presentations, e.g. via spiders [86], spin networks [74] or plane graphs [108].

**Remark 1.5.** Digon removal is used to *evaluate closed webs* in $\mathcal{W}^\eta_k$, i.e. diagrams with boundary labels equal to $k$ only can be simplified to a $\mathbb{Z}[v^{\pm 1}]$-multiple of identity diagrams.

**Remark 1.6.** The category $\mathcal{W}^\eta_2$ (allowing labels 1 and 2) is a $\mathfrak{gl}_2$-analogue of the usual Temperley–Lieb category attached to $\mathfrak{sl}_2$ (where 2 equals the (empty) unit object).

We connect now $\mathcal{W}^\eta$ with *(quantised) fundamental representations* or *(quantised) exterior powers* $\bigwedge^d W$, $d \geq 1$, of the $U_v(\mathfrak{gl}_k)$-modules $W = V_\pm$. The latter is zero if $d > k$ and otherwise defined as the simultaneous $(-\gamma^{-1})$-eigenspace inside $W^{\otimes d}$ for the action of the braid group generators via (1.6). It has the expected explicit basis, namely

$$
(1.12)\quad e_{\mathbf{i}} := e_{i_1} \wedge \cdots \wedge e_{i_d} := \sum_{w \in S_d} (-\gamma)^{-\ell(w)} e_{w(i_1)} \otimes \cdots \otimes e_{w(i_d)}, \quad (k \geq i_1 > \cdots > i_d \geq 1),
$$

indexed by $d$-tuples $\mathbf{i}$. For $\eta = \pm 1$ let $\mathcal{F}und^\eta_k$ be the monoidal category generated by all non-trivial exterior powers of $W$, i.e. objects are tensor products of $\bigwedge^d W$, $1 \leq d \leq k$ inclusively the empty product as monoidal unit. Important morphisms are *q-wedging* and *q-shuffling*:

$$
(1.13)\quad \text{`}q\text{-wedging'}\quad \wedge^{a+b}_{a,b} : \bigwedge^a W \otimes \bigwedge^b W \leftrightarrows \bigwedge^{a+b} W : \vee^{a,b}_{a+b} \quad \text{`}q\text{-shuffling'}.
$$

For example, with $\nu := \delta_{\eta,-1} v^{k-1}$, $e := e_{(1,\ldots,k)}$ and $e(s)$ the same tuple but with $s$ omitted,

$$
(1.14)\quad
\begin{aligned}
\wedge^k_{1,k-1}(e_s \otimes e(j)) &= \delta_{s,j}(-\gamma)^{s-k} \nu e, & \vee^{1,k-1}_k(e) &= \sum_{s=1}^k (-\gamma)^{s-1} \nu^{-1} e_s \otimes e(s), \\
\wedge^k_{k-1,1}(e(j) \otimes e_s) &= \delta_{s,j}(-\gamma)^{1-s} \nu e, & \vee^{k-1,1}_k(e) &= \sum_{s=1}^k (-\gamma)^{k-s} \nu^{-1} e(s) \otimes e_s.
\end{aligned}
$$

We have (for any $k$) the *smoothing relation* $\vee_2^{1,1} \circ \wedge_{1,1}^2 + \gamma \, \text{id} = H_1$, see (1.6). This directly implies with quantised Schur-Weyl duality the first part of the following (where $\eta \in \{\pm 1\}$):

**Proposition 1.7.** *There is a dense full monoidal functor $\Phi_\eta : \mathcal{W}^\eta \to \mathcal{F}und_k^\eta$ which sends a generating object $d$ to $\bigwedge^d W$ and a generating web from (1.9) to the corresponding q-wedging respectively q-shuffling. It induces a monoidal equivalence $\mathcal{W}_k^\eta \simeq \mathcal{F}und_k^\eta$.*

This result provides a purely diagrammatic description of $\mathcal{F}und_k^\eta$. It implies in particular that the non-symmetric braiding morphisms can be expressed in terms of webs, e.g.

$$(1.15) \quad \beta_{a,b} : \bigwedge^a W \otimes \bigwedge^b W \to \bigwedge^b W \otimes \bigwedge^a W, \quad \begin{cases} \sum_{r=0}^a \gamma^{a-r} \, \boxminus \, (a,r,b) & \text{if } a \leq b, \\ \sum_{r=0}^b \gamma^{b-r} \, \boxminus \, (a,r,b) & \text{if } a \geq b. \end{cases}$$

Proposition 1.7 is a reformulation of results from [28]. The authors work in fact with a larger pivotal version, where in the target category one also includes the duals and in the source additionally incorporates flow-lines on webs. From the perspective of tangle invariants it suffices to work with $\mathcal{F}und_k^\eta$ by a clever trick. Namely, we can copy the *RT* construction above, but replace $W^*$ with $\bigwedge^{k-1} W$, the trivial representation with $\bigwedge^k W$ and the cup and cap by the morphisms (1.14). This then provides a monoidal functor $\text{wRT}^\eta$ from oriented tangles to $\mathcal{F}und_k^\eta$. An advantage of this construction is that it stays completely inside $\mathcal{F}und_k^\eta$ and avoids taking duals. This simplifies the situation from a categorification point of view, see Remark 2.51. The invariant of an oriented link $L$ is an endomorphism $f$ of a tensor product of $k$-th exterior powers. Evaluation at the tensor product of the top degree wedges $e$ gives an element in $\mathbb{C}(v)$ which agrees with $J(L)$ in case of $\mathcal{F}und_k^-$ for $k = 2$.

**Example 1.8.** To compute $f = f_{\text{Hopf}}$ for the Hopf link we first translate nested cups (and similarly caps) into webs. Each cup gives rise to copies of $\bigwedge^k W$ depending on the size (indicated by dotted lines). Reading through the webs, see Overview 1, defines the morphism

$$(1.16) \quad \rightsquigarrow \quad f_{\text{Hopf}} : \bigwedge^2 W \otimes \bigwedge^2 W \xrightarrow{\cup} W^{\otimes 4} \xrightarrow{H_1^{-2}} W^{\otimes 4} \xrightarrow{\cap} \bigwedge^2 W \otimes \bigwedge^2 W$$
$$e \otimes e \longmapsto v^3[4] e \otimes e.$$

The construction of this invariant and the more general HOMFLY-PT polynomial via webs and exterior powers is due to [108], [114]. As common in the literature, they work with quantum $\mathfrak{sl}_k$ which produces the same invariant as $\text{wRT}^-$. We prefer to use $\mathfrak{gl}_k$, mainly to make categorifications functorial, see e.g. [43]. In addition, the weight combinatorics and branching rules are much easier, but most importantly, skew Howe duality holds.

**Skew Howe duality.** The crucial observation behind Proposition 1.7 is that a quantum version, established in [28], [90], of a classical tool from invariant theory, namely *skew Howe duality*, can be used to describe all morphisms in $\mathcal{F}und_k^\eta$:

**Proposition 1.9** (q-skew Howe duality). *There is an isomorphism of $U_v(\mathfrak{gl}_k)$-modules*

$$(1.17) \quad \bigwedge^\bullet (W \otimes \mathbb{C}(v)^m) \cong \bigoplus_{\mathbf{d} \in \mathbb{N}_0^m} \bigwedge^{\mathbf{d}} W \quad \text{with} \quad \bigwedge^{\mathbf{d}} W := \bigwedge^{d_1} W \otimes \cdots \otimes \bigwedge^{d_m} W.$$

By turning $\mathbb{C}(v)^m$ into a $U_v(\mathfrak{gl}_m)$-module, one gets commuting mutually centralising actions

(1.18) $\quad U_v(\mathfrak{gl}_k) \quad \curvearrowright \quad X_\pm := \bigwedge^\bullet(V_{\pm,\mathfrak{gl}_k} \otimes V_{\pm,\mathfrak{gl}_m}) \quad \curvearrowleft \quad U_v(\mathfrak{gl}_m)^{\mathrm{op}}$.

The bimodules $X := X_\pm$ inherit some nice symmetry. Namely, the weight spaces for the $U_v(\mathfrak{gl}_m)$-action are direct summands for the $U_v(\mathfrak{gl}_k)$-action and vice versa. As labelling sets we can use their classical weight, i.e. $m$-tuples (respectively $k$-tuples) **c** of integers. Such tuples in fact also index the summands of $X$ from (1.17) and, indeed, there is an isomorphism of vectors spaces $\mathbf{d}^\pm X_\mathbf{c} \cong {}_\mathbf{d}X^{\mathbf{c}^\pm}$. Here, the indices at the top encode the summands and at the bottom the weight space. The left and right position refers to $U_v(\mathfrak{gl}_k)$ and $U_v(\mathfrak{gl}_m)$ respectively, and $\mathbf{d}^\pm$ just means we reverse the tuple $\mathbf{d}$ if the module is $V_-$.

**Example 1.10.** Let $k = 3$ and $m = 2$. Then $\bigwedge^3(V_{\pm,\mathfrak{gl}_3} \otimes V_{\pm,\mathfrak{gl}_2})$ is, as $U_v(\mathfrak{gl}_k)$-module, isomorphic to the direct sum of modules. Each summand becomes a weight space for $U_v(\mathfrak{gl}_m)$ with the action of the generators $E = E_1$, $F = F_1$ indicated via webs:

(1.19)

The labels on the webs encode the $\mathfrak{gl}_m$-weights.

**Remark 1.11.** The $q$-skew Howe duality describes naturally the action of $E_i$ and $F_i$ after projection $1_\mathbf{d}$ onto a weight space. These projections can be encoded conceptually by passing from $U_v(\mathfrak{gl}_m)$ to Lusztig's idempotent version $\dot{U}_v(\mathfrak{gl}_m)$ of $U_v(\mathfrak{gl}_m)$ [95], where idempotent generators $\dot{1}_\mathbf{d}$ are added such that weight modules of $U_v(\mathfrak{gl}_m)$ correspond to modules for $\dot{U}_v(\mathfrak{gl}_m)$. The *fundamental problem* in invariant theory of determining the kernels of the actions is easy in terms of $\dot{U}_v(\mathfrak{gl}_n)$, $n \in \{k, m\}$. By [28], the kernels are the ideals $I_k$ and $I_m$ respectively, generated by all $\dot{1}_\mathbf{d}$, where $\mathbf{d}$ falls outside the respective weight support of $X$.

*Altogether, $\mathcal{F}und_k^\eta$ including its action, braiding and corresponding WRT-tangle invariants is completely controlled by actions of (Lusztig's idempotent version) of quantum groups.*

**Remark 1.12.** There exist variants of $q$-skew Howe dualities, e.g. versions for *i*) symmetric powers [122], *ii*) general linear Lie *super*algebras [116], [140], *iii*) orthogonal and symplectic Lie algebras [128] (replacing $\mathbb{H}_n$ by some Brauer algebra), or *iv*) quantum symmetric pairs (replacing $\mathbb{H}_n$ with a Hecke algebras of Coxeter types $BCD$) [41]. In iii) the dual partner is only a quantum symmetric pair for a fixed point Lie algebra of Langlands dual type inside $\mathfrak{gl}_{2m}$, see [128]; and iv) involves two quantum symmetric pairs for the fixed point Lie algebras $\mathfrak{gl}_k \oplus \mathfrak{gl}_k \subset \mathfrak{gl}_{2k}$ and $\mathfrak{gl}_m \oplus \mathfrak{gl}_m \subset \mathfrak{gl}_{2m}$, [41]. This version sits nicely between the ones from Proposition 1.9 for $(\mathfrak{gl}_k, \mathfrak{gl}_{2m})$ and $(\mathfrak{gl}_{2k}, \mathfrak{gl}_m)$ via restriction/inclusion. It is connected via Hecke algebras of types $BC$ with invariants of knots in an annulus or a disc with a puncture, [57], [58]. A disc with an order two orbifold point can be treated using type $D$ following [5].

**Q1**: *Can Hecke algebras of complex reflection groups treat orbifold points of any order?*

## 2. Four representation theoretic approaches to categorifications

We now sketch categorifications of link and tangle invariants related to the four different views on the Jones polynomial.

### 2.1. Ad I: Khovanov homology

The first categorification of link invariants is given in the work of Khovanov [80] and assigns to an oriented link $L$ a complex $\mathrm{Kh}(L)$ of finite-dimensional graded $\mathbb{C}$-vector spaces. It realises the Jones polynomial $\mathrm{J}(L)$ as the graded Euler characteristic $\chi(\mathrm{Kh}(L))$ of $\mathrm{Kh}(L)$. Thus, $\mathrm{Kh}(L)$ relates to the Jones polynomial of $L$ as a topological space relates to its Betti numbers. Stipulated by the Kauffman bracket, $\mathrm{Kh}$ assigns to the unknot a graded vector space $A$, viewed as complex concentrated in homological degree zero, with Poincaré polynomial $v + v^{-1} = [2]$. Each additional crossing produces a complex one step longer. To make the assignment well-defined one has to work in the homotopy category $K^b(\mathbb{C}\text{-mod}^{\mathbb{Z}\langle v \rangle})$ of the category $\mathbb{C}\text{-mod}^{\mathbb{Z}\langle v \rangle}$ of finite-dimensional $\mathbb{Z}$-graded vector spaces (with grading shift $v$).

Then the *Khovanov invariant* is an assignment

$$\mathrm{Kh}: \{\text{oriented links in } \mathbb{R}^3 \text{ up to isotopy}\} \to K^b(\mathbb{C}\text{-mod}^{\mathbb{Z}\langle v \rangle}) \quad \text{with } \chi(\mathrm{Kh}(L)) = \mathrm{J}(L).$$

Its cohomology, the *Khovanov (co)homology*, and $\mathrm{P}_{\mathrm{Kh}}(L) := \sum_{d,j \in \mathbb{Z}} \dim \mathrm{H}^j(\mathrm{Kh}(L))_d \, t^j v^d \in \mathbb{Z}[v^{\pm 1}, t^{\pm 1}]$, the *Khovanov polynomial*, are invariants as well. The definition of $\mathrm{Kh}$ relies on a *categorified Kauffman bracket* $D \mapsto [\![D]\!]_{\mathrm{cat}}$ with values in $K^b(\mathbb{C}\text{-mod}^{\mathbb{Z}\langle v \rangle})$ whose Euler characteristic is the Kauffman bracket from Section 1. I. This bracket is characterised by

i) the *multiplicativity property* $[\![D_1 \sqcup D_2]\!]_{\mathrm{cat}} = [\![D_1]\!]_{\mathrm{cat}} \otimes_{\mathbb{C}} [\![D_2]\!]_{\mathrm{cat}}$,

ii) the *normalisation* $[\![\bigcirc]\!]_{\mathrm{cat}} = A$, and

iii) the *local smoothing complex* $[\![\times]\!]_{\mathrm{cat}} = [\![\stackrel{\smile}{\frown}]\!]_{\mathrm{cat}} \stackrel{\delta}{\leftarrow} [\![\,|\,|\,]\!]_{\mathrm{cat}}[1]\langle 1 \rangle$.

Local smoothing means that the bracket of a diagram involving a crossing can be expressed as the total complex of a 2-term complex with entries in $K^b(\mathbb{C}\text{-mod}^{\mathbb{Z}\langle v \rangle})$ given respectively by the bracket of the first and the second smoothing (with the shift suggested by $-v = (-1)v^1$ in (1.1)). Since this decreases the number of crossings, by induction one may reduce to the case of no crossing (i.e. circles only), where the functor is specified by *i*) and *ii*). For the construction of the differential $\delta$, Khovanov identifies $A\langle 1 \rangle$ with $\mathbb{C}[x]/(x^2) = H(\mathbb{CP}^1)$ (with $x$ in degree 2), which has additionally a Frobenius algebra structure. The (co)multiplication provide maps $A \otimes A \stackrel{m^*}{\leftarrow} A\langle 1 \rangle$, $A \stackrel{m}{\leftarrow} A \otimes A\langle 1 \rangle$. Applied locally (with appropriate sign rules) they define a map $\delta$ which is then a differential due to the Frobenius algebra properties and sign choices. As for the Jones polynomial one obtains from the bracket a link invariant after incorporating appropriate shifts, i.e. $\text{Kh}(D) = [\![D]\!][n_-]\langle n_+ - 2n_-\rangle$.

Khovanov homology is, as expected, a stronger invariant than the Jones polynomial. Even more striking, Khovanov [78] and Jacobsson [68] could prove that a surface bounded by two links induces a chain map between the Khovanov complexes defining an invariant of the surface, up to signs. The sign issue is fixed in various ways, in [19] via foams, in [33] via surfaces with disorientation lines and in [43] via a sign adaption of Khovanov's construction. The latter, see also [13], provides an *explicit* sign adaption of the involved differential.

**Theorem 2.1.** *The sign adjusted construction of the Khovanov invariant defines a* functor

(2.1)
$$\text{Kh}_{\text{sgn}} : \{\textit{oriented links in } \mathbb{R}^3\} \to K^b(\mathbb{C}\text{-mod}^{\mathbb{Z}(v)}) \quad \textit{with homology} \quad \text{P}_{\text{Kh}_{\text{sgn}}} = \text{P}_{\text{Kh}}.$$

This *functoriality* is crucial for topological applications, e.g. to prove Milnor's conjecture on slice genus of torus knots [118] or unknot detection by Khovanov homology [85].

**Remark 2.2.** The (categorified) Kauffman bracket works well for links. For tangles, an additional direction of composition has to be reflected in the target category of a possible invariant. Instead of working with the category of vector spaces, one has to pass to e.g. categories of bimodules over (generalised) Khovanov arc algebras [79], [134], [30], operads and canopolis [10], or various topological incarnations related to foam categories. An analogue, although not very practical, of the Kauffman bracket for $\mathfrak{gl}_k$, $k > 2$, can be given via (1.15).

**Remark 2.3.** In practice one often considers all complete smoothings at once and arranges their values as vertices in the famous cube of resolutions [80], [11], with the differential on the edges. However, the interpretation of (1.1) in terms of a 2-term complex was chosen to highlight the important role played by such complexes in algebraic(-geometric) categorifications. They do not only appear in crucial definitions (like coherent sheaves or other Serre quotient categories), but also provide technical toolkits for categorical actions for instance in form of spherical twists, spherical functors or Rickard complexes.

**Remark 2.4.** Although *odd Khovanov homology* and *Lee homology* are often called variants of Kh, they are rather different theories from our point of view. Instead of $\mathfrak{gl}_2$, they are connected with the Lie (super) algebras $\mathfrak{osp}(1|2)$ [50] and $\mathfrak{gl}_1 \times \mathfrak{gl}_1$ [123], respectively.

**Q2** : *Which surfaces are distinguished by Khovanov homology? Does it provide a* 4-*TQFT?*

### 2.2. Ad II: Triply graded link homology

A categorification of link invariants using braid closures and traces is the *triply graded Khovanov–Rozansky homology*. It was originally constructed using matrix factorisations [82] and then reinterpreted [77] in representation theoretic terms via Soergel bimodules. To sketch the construction, we view $\beta \in \text{Br}_n$ as a special tangle with $n$ bottom and $n$ top points or as a "map" from $n$ inputs to $n$ outputs. We associate variables $x_1, \ldots, x_n$ to the inputs and consider the category $R_n$-mod of modules over the polynomial ring $R_n =$ $\mathbb{C}[x_1, \ldots, x_n]$. To $\beta$ we assign a certain (complex of) $R_n$-bimodule(s) $X(\beta)$ which defines a "map" $X(\beta) \otimes_{R_n} \_ : R_n\text{-mod} \to R_n\text{-mod}$. Taking the closure $\hat{\beta}$ of $\beta$ connects or identifies the points at the bottom with those at the top, see Example 1.1. Categorically this corresponds to identifying the left with the right action of $R_n$ on $X(\beta)$. Algebraically one takes (derived) coinvariants, i.e. Hochschild homology of $X(\beta)$. This Hochschild homology is a bigraded vector space with gradings coming from the Hochschild and homological grading. It is even triply graded if one works with graded $R_n$-modules.

To be more rigorous, consider $R := R_n$ as a graded ring with $\deg(x_i) = 2$ and let $S_n$ act on $R$ by permuting the variables. Given any subset $I \subset S_n$ of simple transpositions, let $R^I = R^{W_I} \subset R$ be the ring of invariants under the action of $I$ or equivalently under the parabolic subgroup $W_I$ generated by $I$ inside $W = S_n$.

**Example 2.5.** *Obviously, $R^W$ is the ring of symmetric polynomials and $R^\emptyset = R$. In case $I = \{s_i\}$ we obtain $R^{s_i} := R^{\{s_i\}} = \mathbb{C}[x_1, \ldots, x_{i-1}, x_i + x_{i+1}, x_i x_{i+1}, x_{i+2}, \ldots, x_n]$.*

To any word $\ddot{w} = s_{i_1} s_{i_2} \ldots s_{i_r}$ in simple transpositions from $S_n$, there is an associated *Soergel bimodule* $\mathrm{BS}(\ddot{w}) = R \otimes_{R^{s_{i_1}}} R \otimes_{R^{s_{i_2}}} \cdots \otimes_{R^{s_{i_r}}} R \langle -\ell(\ddot{w}) \rangle$ which is the *Bott–Samelson bimodule for $\ddot{w}$*, see Remark 2.18. In particular, $\mathrm{BS}(s_i) = R \otimes_{R^{s_i}} R \langle -1 \rangle$ and $\mathrm{BS}(\emptyset) = R$.

The *category of Soergel bimodules* $\mathcal{SBim}_n$ [132], [131] is defined as the Karoubian closure of the additive category generated by Bott–Samelson bimodules and its grading shifts (inside the category of all graded bimodules with degree zero maps). It is an additive category and closed under $\_ \otimes_R \_$, i.e. it is monoidal with unit $R$. Next, the generating *Rouquier complexes* associated with $\beta_i$ and $\beta_i^{-1} \in \mathrm{Br}_n$ are

$$X(\beta_i): \quad \left( R\langle 1 \rangle \longrightarrow \boxed{\mathrm{BS}(s_i)} \right), \quad \text{and} \quad X(\beta_i^{-1}): \quad \left( \boxed{\mathrm{BS}(s_i)} \longrightarrow R\langle -1 \rangle \right)$$

with differentials given by $1 \mapsto (x_i - x_{i+1}) \otimes 1 + 1 \otimes (x_i - x_{i+1})$ and $1 \otimes 1 \mapsto 1$ respectively.

To an element $\beta \in Br_n$ written as a word $\ddot{\beta} = \beta_{i_1}^{\epsilon_1} \cdots \beta_{i_r}^{\epsilon_r}$ in the $\beta_i^{\pm 1}$, Rouquier [124] attaches the corresponding tensor product $X(\ddot{\beta}) := X(\beta_{i_1}^{\epsilon_1}) \otimes_R \cdots \otimes_R X(\beta_{i_r}^{\epsilon_r})$ (with the convention that the identity braid in $Br_n$ is sent to $R$) in $K^b(\mathcal{SBim}_n)$. He then proves the following important result which allows one to use the notation $X(\beta)$.

**Theorem 2.6.** *If two words $\ddot{\beta}$ and $\ddot{\beta}'$ represent the same element in $Br_n$, then the Rouquier complexes $X(\ddot{\beta})$ and $X(\ddot{\beta}')$ are canonically isomorphic in $K^b(\mathcal{SBim}_n)$.*

**Remark 2.7.** Rouquier [124] in fact constructed a *genuine braid group action* on $K^b(\mathcal{SBim}_n)$ by these Rouquier complexes. Explicit rigidity maps (even over $\mathbb{Z}$) were determined in [47].

To categorify braid closures and traces, consider the Hochschild homology functor

$$\mathrm{HH}(\_) := \bigoplus_{i \in \mathbb{N}_0} \mathrm{HH}_i(R, \_) := \bigoplus_{i \in \mathbb{N}_0} \mathrm{Tor}_i(R, \_)$$

from the category of $\mathbb{Z}$-graded $R$-bimodules to the category of $(\mathbb{N}_0 \times \mathbb{Z})$-graded (viewed as $\mathbb{Z} \times \mathbb{Z}$-graded) vector spaces. For a complex $C$ of finitely generated graded $R$-bimodules, let $\mathrm{HH}(C)$ denote the complex of bigraded abelian groups obtained by applying the functor $\mathrm{HH}$ to the components and differentials of $C$. Set $\mathbb{HHH}(C) := \mathrm{H}^\bullet(\mathrm{HH}(C))$. This is an object in the category $\mathbb{C}\text{-mod}^{\mathbb{Z}(t,v,h)}$ of triply graded vector spaces with $t$, $v$ and $h$ referring to the homological, internal and Hochschild degree respectively (with shift functors $[\_]$, $\langle\_\rangle$, $\{\_\}$). Its *(3-parameter) Poincaré series* is a Laurent series in $v$ with coefficients in $\mathbb{Z}[h^{\pm 1}, t^{\pm 1}]$:

$$(2.2) \quad \mathrm{P}(\mathbb{HHH}(C)) := \sum_{d,i,j \in \mathbb{Z}} \dim(\mathrm{H}^j(\mathrm{HH}_i(C))_d) \, t^j h^i v^d \ \in \ \mathbb{Z}[h^{\pm 1}, t^{\pm 1}]((v)).$$

We have to work here with *Laurent series* in $v$, indicated by $((v))$, since $R$ is infinite-dimensional, but the expression makes sense since the components are finite in each fixed triple degree. Evaluating $t = -1$ gives the *graded Euler characteristic* $\chi(\mathbb{HHH}(C)) \in \mathbb{Z}[h^{\pm}]((v))$. Khovanov showed in **[77]** that (2.2) gives, up to some rescaling, an invariant of oriented links (here $\epsilon(\ddot{\beta})$ denotes the sum of the exponents of the $\beta_i$ appearing in $\ddot{\beta}$):

**Theorem 2.8.** *For a braid word $\ddot{\beta}$ in $Br_n$, the normalised Poincaré series*

$$(2.3) \quad \mathrm{KR}(\hat{\beta}) := (th)^{\frac{1}{2}(\epsilon(\ddot{\beta})-n)} v^{\epsilon(\ddot{\beta})} \, \mathrm{P}(\mathbb{HHH}(X(\ddot{\beta})))$$

*only depends on the braid closure $\hat{\beta}$. Thus there is a well-defined assignment*

$$(2.4) \quad \mathrm{KR} : \{\textit{oriented links in } \mathbb{R}^3 \textit{ up to isotopy}\} \to \mathbb{Z}[h^{\pm\frac{1}{2}}, t^{\pm\frac{1}{2}}]((v)), \quad \hat{\beta} \mapsto \mathrm{KR}(\hat{\beta}).$$

*The invariant* $\mathrm{KR}(\hat{\beta})$ *is called the* triply graded Khovanov–Rozansky homology *of $\hat{\beta}$.*

**Example 2.9.** We calculate $\mathrm{KR}(\bigcirc)$, i.e. $\beta$ is for $n = 1$ the identity braid in $Br_n$ with $X(\beta) = R_n = R = \mathbb{C}[x_1]$. The Hochschild homology can easily be computed as $\mathrm{HH}_0(R) = R$, $\mathrm{HH}_1(R) = R\langle 2 \rangle$, $\mathrm{HH}_{\geq 2}(R) = \{0\}$ from the Koszul resolution $\mathbb{C}[y] \otimes \mathbb{C}[y'] \stackrel{(y-y')\cdot}{\to} \mathbb{C}[y] \otimes \mathbb{C}[y']$ of $R = \mathbb{C}[x_1]$ (with $y, y' \mapsto x_1$). Since the Poincaré series of $R$ equals $\mathrm{P}(R) = \frac{1}{1-v^2}$, we obtain

$$(2.5) \quad \mathrm{KR}(\bigcirc) = t^{-\frac{1}{2}} h^{-\frac{1}{2}} \mathrm{P}(\mathbb{HHH}(R)) = t^{-\frac{1}{2}} h^{-\frac{1}{2}} \frac{1 + hv^2}{1 - v^2} \in \mathbb{Z}[h^{\pm}, t^{\pm}]((v)).$$

For general $n$, $\mathrm{HH}(R_n) \cong R_n \otimes \bigwedge^{\bullet}(\xi_1, \ldots, \xi_n)$, where each $\xi_i$ is of $v$-degree 2 and $h$-degree 1, and $\mathbb{HHH}(R_n) = \mathrm{HH}(R_n)$. As expected, the identity braid gives then, with (2.3), $\mathrm{KR}(\bigcirc)^n$.

Crucial for the proof of Theorem 2.8 are isomorphisms $\mathbb{HHH}(X(\alpha) \otimes_R X(\beta_n)) \cong \mathbb{HHH}(X(\alpha))\langle -1 \rangle$, $\mathbb{HHH}(X(\alpha) \otimes X(\beta_n^{-1})) \cong \mathbb{HHH}(X(\alpha))\langle 1 \rangle[1]\{1\}$ and the trace property $\mathbb{HHH}(X(\alpha) \otimes_R X(\alpha')) \cong \mathbb{HHH}(X(\alpha') \otimes_R X(\alpha))$ for $\alpha, \alpha' \in Br_n \subset Br_{n+1}$. The latter follows directly from the canonical isomorphism $\mathrm{HH}(M \otimes_R N) = \mathrm{HH}(N \otimes_R M)$ noting that for Soergel bimodules one does not need to derive the tensor product, since they are by standard invariant theory free as one-sided modules. The formulas imply that $\hat{\beta} \to \chi(\mathbb{HHH}(X(\beta)))$ factors through a $\mathbb{Z}[v^{\pm 1}]$-linear trace function $\tau : \coprod_{n \geq 1} \mathbb{H}_n \to \mathbb{Z}[h^{\pm 1}]((v))$ such that

$$(2.6) \quad \tau(1) = \frac{1 + hv^2}{1 - v^2}, \quad \tau(xH_n^{\pm 1}) = z_{\pm}\tau(x) \ \forall x \in \mathbb{H}_n \text{ with } z_+ = v^{-1} \text{ and } z_- = -hv,$$

where $H_i$ is the image of $\beta_i$ in $\mathbb{H}_n$. With the normalisation (2.3), $\tau$ becomes a Markov trace. By introducing a homological $\frac{1}{2}\mathbb{Z}$-grading, even $\mathbb{HHH}$ can be turned into an invariant, see **[124]**.

**Remark 2.10.** Since trace functions $\tau$ on $\coprod_{n \geq 1} \mathbb{H}_n$ are classified by the pair $(\tau(1), z_+)$ as in (2.6), one can identify $\tau$ from (2.6), up to normalisation of $\tau(1)$ with the (Jones–)Ocneanu trace **[70]**. In **[87]** it is proved that for any finitely generated Coxeter group, with the more

14    ICM 2022

general definition of Soergel bimodules and Hecke algebra from Theorem 2.19, the Euler characteristic of the KR-homology provides a Markov trace on the Hecke algebra.

**Remark 2.11.** To get nicer formulas we make a change of variables by setting $\mathbf{a} = v(ht)^{\frac{1}{2}}$. Then $\mathrm{KR}(\hat{\beta}) \in \mathbb{Z}[\mathbf{a}^{\pm 1}, t^{\pm 1}]((v))$. For example, $\mathrm{KR}(\bigcirc) = v \frac{\mathbf{a}^{-1} + \mathbf{a}t^{-1}}{1 - v^2} = \frac{\mathbf{a}^{-1} + \mathbf{a}t^{-1}}{v^{-1} - v}$. Setting $t = -1$ gives $\frac{\mathbf{a} - \mathbf{a}^{-1}}{v - v^{-1}}$ and the characterising skein relation (1.2) of the HOMFLY-PT polynomial holds. Setting $\mathbf{a} = v^2$ we obtain the Jones polynomial, e.g. $v + v^{-1}$ here and $v^3[4]$ in Example 2.13.

**Theorem 2.12.** $\mathrm{KR}(\hat{\beta}) \in \mathbb{Z}[\mathbf{a}^{\pm 1}, t^{\pm 1}]((v))$ *specialises for* $t = -1$ *to the* HOMFLY-PT *polynomial.*

**Example 2.13.** For $L$ the Hopf link from (1.1) we get $\chi(\mathrm{KR}(L)) = v^2 \tau(H_1^2)$. By (1.7), $\tau(H_1^2) = \tau(1) + (v^{-1} - v)\tau(H_1) = \sigma^2 + (v^{-1} - v)v^{-1}\sigma$ with $\sigma := \tau(1)$.

**Remark 2.14.** The appearance of $\mathbb{H}_n$ here is not surprising, because Soergel originally invented his bimodules to understand the Kazhdan–Lusztig basis in the Hecke algebra $\mathbb{H}_n$. To formulate this more precisely, let $\mathrm{K}_0^{\oplus}(\mathcal{SBim}_n)$ be the split Grothendieck ring of the additive category of Soergel bimodules. That is the free $\mathbb{Z}[v^{\pm 1}]$-module generated by isomorphism classes $[M]$ of objects $M$ in $\mathcal{SBim}_n$ modulo relations $[M \oplus N] = [M] + [N]$, $v[M] = [M\langle 1 \rangle]$ and multiplication $[M][N] = [M \otimes_R N]$. In [132], Soergel proved an influential categorification theorem which is crucial for all representation theoretic constructions of categorified link invariants: there is an isomorphism of $\mathbb{Z}[v^{\pm 1}]$-algebras

$$(2.7) \qquad \Upsilon_n : \mathbb{H}_n \longrightarrow \mathrm{K}_0^{\oplus}(\mathcal{SBim}_n), \quad H_i + v \longmapsto [B(i)].$$

which moreover identifies the Kazhdan–Lusztig basis with classes of indecomposables bimodules. Hereby $H_i$ corresponds to a virtual object only. This can be fixed by identifying $\mathrm{K}_0^{\oplus}(\mathcal{SBim}_n)$ with the Grothendieck group of the triangulated category $K^b(\mathcal{SBim}_n)$, since then $[X(\beta_i)] = \Upsilon(H_i)$. This shows that Rouquier's braid group action, despite its faithfulness [83], is honestly a categorical Hecke algebra action which also descends to a Hecke algebra action on the Grothendieck group. The relations (1.1), (1.2), (1.7) indicate that the presented invariants should be rather connected with the Hecke algebra instead of the braid group.

In contrast to Kh, computing KR is usually hard, although the resulting values might be more conceptual and expressible using generating series. Important progress was however made recently for the *torus links* $t_{(p,q)}$ which are the closure of $(\beta_{p-1} \cdots \beta_2 \beta_1)^q$ (with the Hopf link as special case $(p, q) = (2, 2)$). An important first step is done in [66] with the observation that $\mathrm{KR}(t_{(n,q)})$ stabilises for $q \to \infty$ to a limit isomorphic to $\mathbb{C}[u_1, \ldots, u_n] \otimes \bigwedge^{\bullet}[\xi_1, \ldots \xi_n]$ with $u_i$ in $h$-degree zero and $\xi_i$ in $h$-degree 1 (cf. $\mathrm{HH}(R_n)$ in Example 2.9). This limit is identified in [66] with the derived endomorphism ring of a certain categorified Young idempotent in the Hecke algebra $\mathbb{H}_n$. This idempotent provides a bridge to categorified coloured *RT*-invariants, Remark 2.54 and Conjecture 2.58, since it acts on $(V_{-,\mathfrak{gl}_2})^{\otimes n}$ as a projector, the $V_{-,\mathfrak{gl}_2}$-version of (2.21).

In [45], $\mathrm{KR}(t_{(p,q)})$ is determined via a beautiful recursive formula in case $p = q$, and extended to general $(p, q)$ in [67]. They both use categorifications of idempotents in $\mathbb{H}_n$ which are interesting tools on their own, e.g. for developing a categorified representation theory of $\mathbb{H}_n$. For general links, computing KR seems still to be out of reach. Instead of studying the invariant via its original definitions [77], [82], alternative constructions were proposed, e.g. the following involving Hilbert schemes and Cherednik algebras, inspired in particular by the underlying combinatorics of symmetric functions and Macdonald polynomials.

**Remark 2.15.** The approach of [110], [109] starts by viewing a torus link $L = t_{(p,q)}$ as an algebraic link, i.e. as the intersection of a planar curve $C := C_{p,q} \subset \mathbb{C}^2$ (defined by the polynomial $f = x^p - y^q \in \mathbb{C}[x,y]$) with a sufficiently large sphere around the origin in $\mathbb{C}^2$. Attached to $C$ is the *Hilbert scheme* $C^{[r]}$ of $r$ points on $C$ which, as a set, is given by all ideals $I \subset \mathbb{C}[x,y]$ of codimension $r$ containing $f$. In [110] it is proved for coprime $p, q$ that the Euler characteristic of $\mathrm{KR}(t_{(p,q)})$, i.e. the HOMFLY-PT polynomial, equals up to a normalisation the Euler characteristic of the disjoint union of all *nested Hilbert schemes*

$$C^{[d,d+i]} := \{(I, J) \mid I \cdot (x, y) \subset J \subset I\} \subset C^{[d]} \times C^{[d+i]}$$

with $d$ and $i$ encoding the $v$- respectively **a**-degree. For a generalisation to algebraic links, see [104]. In [109] it is conjectured that replacing the Euler characteristic with the virtual Poincaré polynomial, see loc. cit. for the definition, provides the triply graded Khovanov homology. For torus links, this is proved in [111]. In general, this conjecture is still open.

**Remark 2.16.** As indicated in the introduction, KR is related to double affine Hecke algebras (DAHAs) and their rational degenerations from [51]. The rational DAHA $\mathrm{H}_c = \mathrm{H}_c(S_n)$ with parameter $c \in \mathbb{C}$ is the quotient of $\mathbb{C}\langle x_i, y_i \mid 1 \leq i \leq n \rangle \rtimes S_n$ modulo

$$[x_i, x_j] = 0 = [y_i, y_j], \quad [x_i, y_j] = c \cdot (i, j), \quad [x_i, y_i] = 1 - c \sum_{j \neq i} (i, j)$$

for any $i \neq j$ with $(i, j) \in S_n$. It is a flat deformation of $\mathrm{H}_0 = \mathcal{D} \rtimes S_n$, where $\mathcal{D}$ is the algebra of differential operators on $\mathbb{C}^n$. If $c = \frac{p}{q}$ with $(p, q) \in \mathbb{Z}^2$ coprime, there is a unique irreducible finite-dimensional $\mathrm{H}_{\frac{p}{q}}$-module $\mathrm{L}_{p,q}$ [15]. When restricting the $\mathrm{H}_c$-action to $S_n$, this module decomposes into direct summands. Let $1_i \mathrm{L}_{p,q}$ be the isotypic component of $\bigwedge^i \mathbb{C}^{n-1}$ using the reflection representation $\mathbb{C}^{n-1}$. The internal grading on $\mathrm{H}_c$ realised by the eigenvalues of the Euler operator $\mathrm{eu} = \sum_{i=1}^n x_i y_i$ (and encoding the difference of the polynomial degree in the $x$'s and the $y$'s) induces a grading on $\mathrm{L}_{p,q}$ and $1_i \mathrm{L}_{p,q}$. In [61], the Poincaré polynomial $P(M)$ of $M := \bigoplus_i 1_i \mathrm{L}_{p,q}$ is identified with the HOMFLY-PT polynomial of $t_{(p,q)}$ up to renormalisation. Here, $i$ contributes to the **a**-degree and eu to the $v$-degree. The identification is achieved by matching known formulas for the HOMFLY-PT polynomial with the character formula for $\mathrm{L}_{p,q}$ from [15]. In [61], a filtration on $M$ is predicted such that $\mathrm{KR}(t_{(p,q)})$ arises as $P(\mathrm{gr}M)$ for the associated graded $\mathrm{gr}M$. This is verified in [111] in terms of a geometric perverse filtration, after realising $1_i \mathrm{L}_{p,q}$ (with the action of the spherical Hecke algebra $1_i \mathrm{H}_{\frac{p}{q}} 1_i$) as $\bigoplus_d H(C^{[d,d+i]})$ (with the action of certain Hecke–Nakajima operators). The comparison of $P(\mathrm{gr}M)$ with $\mathrm{KR}(t_{(p,q)})$ is again done by matching explicit formulas from [45], [67].

**Q3**: *Is there a combinatorial model to compute KR? For which cobordisms is KR functorial?*

**2.2.1. Interlude: Hecke categories.** The quantum $\mathfrak{gl}_k$-invariants and the construction of the fundamental representations (1.12) use heavily the monoidal structure of $\mathcal{R}ep_k$. By (1.7), the action of the braid group on $(V_\pm)^{\otimes n}$ factors through an $\mathbb{H}_n$-action preserving the weight spaces of $(V_\pm)^{\otimes n}$. To get categorified tangle invariants one might therefore categorify these Hecke algebra actions in terms of a monoidal category acting via functors on a category, ideally with an extension to categorified quantum group actions and $q$-skew Howe duality (1.18). To motivate the origin of such actions we go back to the original definition of Hecke algebras arising from split reductive groups G defined over a finite field $\mathbb{F}_q$ with a choice T ⊂ B ⊂ G of a maximal torus and Borel subgroup and the finite group $G(\mathbb{F}_q)$ of $\mathbb{F}_q$-points. Most finite simple groups, in particular of Lie type, arise in this way. For us the case of $G = \mathrm{GL}_n$ suffices with the choice of diagonal matrices inside the upper triangular matrices and their corresponding finite groups $G_q := \mathrm{GL}_n(\mathbb{F}_q) \supset B_q \supset T_q$. The *Weyl group* $W = N_{G_q}(T_q)/T_q$ can be identified with the group $S_n \subset G_q$ of permutation matrices.

The associated *Iwahori–Hecke algebra* $\mathrm{H}_n(q)$ is the vector space $\mathrm{Func}_{B_q \times B_q}(G_q, \mathbb{C})$ of complex valued functions $f$ on $G_q$ invariant under both the left and the right action of $B_q$. i.e. $f(bg) = f(g) = f(gb)$ for all $g \in G_q$, $b \in B_q$, equipped with the convolution product

$$(2.8) \qquad (f \star g)(x) = \frac{1}{|B_q|} \sum_{y \in G_q} f(xy^{-1}) g(y).$$

The indicator functions $h_w$, $w \in W$, for the double cosets $B_q w B_q$ form a basis of $\mathrm{H}_n(q)$ by the *Bruhat decomposition* (or just Gauss elimination) $G_q = \bigsqcup_{w \in W} B_q w B_q$. In this basis, the structure constants of the multiplication are polynomial in $q = |\mathbb{F}_q|$ and thus one can replace $q$ by a generic variable and "treat all $q$ at once". Then $\mathrm{H}_n(q)$ becomes isomorphic to $\mathbb{H}_n$ via $q \mapsto v^{-2}$, $h_{s_i} \mapsto v^{-1} H_i$ after adjoining a square root of $q$.

**Remark 2.17.** The construction allows vast generalisations, e.g. by replacing $\mathbb{F}_q$ by a local field with finite residue field (to get Iwahori–Hecke algebras arising in number theory), by working with topological groups, or with convolution products in homology theories.

The usual Grothendieck function-sheaf correspondence, see e.g. [89], indicates that a categorification is given by a certain category of $(B \times B)$-equivariant sheaves on G. Since $(B \times B)$-equivariant functions on G can be identified with B-equivariant functions on G/B, a categorification might therefore work with B-equivariant sheaves on G/B.

**Overview 2**  The geometric, algebraic and Lie theoretic Hecke categories

| Hecke algebra | ⤳ Hecke category | $\stackrel{\text{Thm 2.19}}{\simeq}$ Soergel bimodules | $\stackrel{\text{Rk 2.24}}{\simeq}$ Projective functors |
|---|---|---|---|
| $\mathbb{H}_n, \mathrm{H}_n(q)$ | $\mathcal{H}_n^{\text{geo}}$ | $\mathcal{SBim}_n$ | $\mathcal{P}_n$ |

For the categorification we use geometry over $\mathbb{C}$ with $G := \mathrm{GL}_n(\mathbb{C})$, $B := \mathrm{B}(\mathbb{C})$, $T := \mathrm{T}(\mathbb{C})$ and the algebraic variety $\mathcal{F} = G/B$ of all full flags $\{F_1 \subset \cdots \subset F_n = \mathbb{C}^n \mid \dim(F_i) = i\}$ of vector subspaces in $\mathbb{C}^n$. The bounded equivariant derived category $\mathcal{D}_B^b(\mathcal{F}, \mathbb{C})$ of sheaves of $\mathbb{C}$-vector spaces [17], is a monoidal category with a convolution product $\star$ [133].

The *geometric Hecke category* $\mathcal{H}_n^{\text{geo}}$ is defined as the full subcategory of $\mathcal{D}_B^b(\mathcal{F}, \mathbb{C})$ generated by the constant sheaves $\underline{\mathbb{C}}_{P_i}$ on $P_i = \overline{B s_i B} = B s_i B \cup B \subset G$ under convolution $\star$, homological shifts [1], finite direct sums and direct summands. Concretely, the objects in $\mathcal{H}_n^{\text{geo}}$ are shifts of objects $\text{BS}^{\text{geo}}(\ddot{w}) = \underline{\mathbb{C}}_{P_{i_1}} \star \cdots \star \underline{\mathbb{C}}_{P_{i_r}} [-r]$ for any word $\ddot{w} = s_{i_1} \ldots s_{i_r}$ in simple transpositions, and finite direct sums and summands of those. The shift functors $[i]$ turn $\mathcal{H}_n^{\text{geo}}$ into a graded category. Note the similarity to $\mathcal{SB}im_n$ with shift functors $\langle i \rangle$.

**Remark 2.18.** The objects $\text{BS}^{\text{geo}}(\ddot{w})$ have a nice alternative description in terms of the *Bott–Samelson varieties* $Z(\ddot{w}) = P_{i_1} \times \cdots \times P_{i_{r+1}}/B^r$, where $y = (y_1, \ldots, y_r) \in B^r$ acts as $y.(p_1, \ldots, p_r) = (p_1 y_1^{-1}, y_1 p_2 y_2^{-1}, \ldots, y_r p_{r+1})$. If $\ddot{w}$ is a reduced expression for $w \in W$, then the multiplication map $\pi : Z(\ddot{w}) \to G/B, (p_1, \ldots, p_r) \mapsto p_1 \cdots p_r$ is known to be a resolution of singularities for the Schubert variety $\overline{BwB}/B$, first studied in the context of compact Lie groups by Bott and Samelson. It is not hard to see that $\text{BS}^{\text{geo}}(\ddot{w}) \cong \pi_* \underline{\mathbb{C}}_{Z(\ddot{w})}$ and that the Bott-Samelson bimodules arise as $T$-equivariant cohomology $H_T(Z(\ddot{w})) \cong \text{BS}(\ddot{w})$.

Soergel's categorification result, Remark 2.14, arises now naturally:

**Theorem 2.19.** *There is an equivalence $\mathcal{H}_n^{\text{geo}} \simeq \mathcal{SB}im_n$ of graded monoidal categories sending* $\text{BS}^{\text{geo}}(\ddot{w})$ *to* $\text{BS}(\ddot{w})$. *In particular,* $K_0^\oplus(\mathcal{H}_n^{\text{geo}}) \cong \mathbb{H}_n$ *as* $\mathbb{Z}[v^{\pm 1}]$-*algebras.*

**Remark 2.20.** Theorem 2.19 can be proved by identifying both, $\mathcal{H}_n^{\text{geo}}$ **[121]** and $\mathcal{SB}im_n$ **[48]**, with the Karoubian closure $\mathcal{DB}im_n$ of a *diagrammatic monoidal category* $\mathcal{DB}im_n'$ invented in **[46]**, **[48]** and proved to be equivalent to the subcategory of Bott–Samelson bimodules. Striking is that this category $\mathcal{DB}im_n'$ has a *presentation* with generators and relations. One of the most prominent applications is the proof of the long outstanding positivity conjecture for the Kazhdan–Lusztig polynomials of an arbitrary Coxeter system and an algebraic proof of the Kazhdan–Lusztig conjectures for reductive complex Lie algebras in **[49]**.

### 2.3. Ad IV: Categorification of the web calculus and its tangle invariant

We reverse the order from Section 1 and pass to categorifications for $\text{wRT}^\pm$ which are further developed than for $\text{RT}$. A categorification of the quantum $\mathfrak{gl}_n$ tangle invariant $\text{wRT}^\pm$ is constructed by Mazorchuk and the author **[135]**, **[105]** and Sussan **[139]**, using highest weight categories of infinite-dimensional representations of the (again, but now in a different role!) general linear Lie algebras $\mathfrak{gl}_N(\mathbb{C})$. It categorifies $\mathcal{F}und^\pm$ and even skew Howe duality: objects $\bigwedge^{\mathbf{d}} V_\pm$ as in (1.17) are realised as Grothendieck groups of categories, and actions and morphisms are lifted to functors with relations realised by specific natural transformations. This construction is part of a major change of perspective in representation theory in recent years. The starting point goes back to Crane and Frenkel **[35]**, who proposed the idea of *Hopf categories* to construct 4-TQFTs based on categorified quantum groups and canonical bases. Categorified quantum groups were then defined in **[81]**, **[125]** as certain 2-categories. We will indicate later how they arise naturally in the context of categorified tangle invariants.

Let $\mathfrak{h} \subset \mathfrak{b} \subset \mathfrak{g} = \mathfrak{g}_N := \mathfrak{gl}_N(\mathbb{C})$ be the Cartan and Borel subalgebra given by all diagonal respectively upper triangular matrices. Equip $\mathfrak{h}^*$ with the standard basis $\delta_1, \ldots, \delta_N$, such that $\delta_i$ picks out the $i$th diagonal matrix entry, and with the symmetric bilinear form $(\delta_i, \delta_j) = \delta_{i,j}$. We identify the lattice $\mathfrak{h}_{\text{int}} := \mathbb{Z}\delta_1 \oplus \cdots \oplus \mathbb{Z}\delta_n$ of integral weights with $\mathbb{Z}^N$ via $\lambda \leftrightarrow (\lambda_1, \ldots, \lambda_N)$, where $\lambda_i = (\lambda + \rho, \delta_i)$ with $\rho = \sum_{j=1}^N (N - j + 1)\delta_j$. The group $S_N$ acts on $\mathbb{Z}^N = \mathfrak{h}_{\text{int}}$ by permuting components and defines the *Bruhat ordering* generated by $\mu < \lambda$ if $\lambda$ differs from $\mu$ by swapping a pair $\mu_i, \mu_j$ such that $\mu_i < \mu_j$ and $i < j$.

We set up now a dictionary between standard basis vectors $\vec{e} \in \bigwedge^{\mathbf{d}} V_{\pm}$ of $\bigwedge^{\mathbf{d}} V_{\pm}$ (for fixed $\pm$) and a subset $\Lambda^{\mathbf{d}} \subset \mathbb{Z}^N$ of $\mathfrak{g}_N$-weights (with $N = \sum_{i=1}^m d_i$). Each tensor product $\vec{e}$ of basis vectors (1.12) is identified with an element $\text{wt}(\vec{e}) \in \{1, 2, \ldots, k\}^N \subset \mathbb{Z}^N$ via

$$(2.9) \qquad \vec{e} = e_{\mathbf{i}}^{(1)} \otimes \cdots \otimes e_{\mathbf{i}}^{(m)} \quad \mapsto \quad \text{wt}(\vec{e}) := (i_1^{(1)}, \ldots, i_{d_1}^{(1)}, i_1^{(2)}, \ldots, i_{d_m}^{(m)}) \in \Lambda^{\mathbf{d}}.$$

Let $\Lambda^{\mathbf{d}}$ be the image. Note that a weight space in $\bigwedge^{\mathbf{d}} V_{\pm}$ corresponds to an $S_N$-orbit $\mathbf{c}$ in $\Lambda^{\mathbf{d}}$.

Now we construct a category $O^{\mathbf{d}}$ whose Grothendieck group has a basis naturally labelled by $\Lambda^{\mathbf{d}}$. For this consider the BGG category $O$ of all finitely generated $\mathfrak{g}$-modules $M$ which are locally finite over $\mathfrak{b}$ and have a weight space decomposition with only integral weights $\lambda \in \mathfrak{h}_{\text{int}}$. This is an abelian finite length category, where the irreducible objects are exactly the irreducible highest weight modules $L(\lambda)$ of highest weight $\lambda \in \mathfrak{h}_{\text{int}}$, i.e. the irreducible quotients of the *Verma modules* $\Delta(\lambda)$ for $\lambda \in \mathfrak{h}_{\text{int}}$. Objects in $O$ which have a $\Delta$-flag, i.e. a finite filtration with subquotients isomorphic to Verma modules, form an exact additive subcategory $O^{\Delta}$ which is closed under direct summands and contains all projective objects. Even more, category $O$ is a *highest weight category*, see e.g. [24], for the set $\mathfrak{h}_{\text{int}} = \mathbb{Z}^N$ viewed as poset with *standard objects* the $\Delta(\lambda)$. Technically this means that the projective cover of $L(\lambda)$ surjects onto $\Delta(\lambda)$, and $\Delta(\lambda)$ surjects onto $L(\lambda)$, and the kernel has a $\Delta$-flag with subquotients some $\Delta(\mu)$ where $\mu > \lambda$, respectively a Jordan–Hölder filtration with subquotients $L(\mu)$'s with $\mu < \lambda$. As a consequence, the canonical maps induce isomorphisms between Grothendieck groups for i) the additive category of projectives, ii) the exact category $O^{\Delta}$, iii) the abelian category $O$ and iv) the triangulated bounded derived category $D^b(O)$:

$$(2.10) \qquad K_0^{\oplus}(\text{Proj}(O)) = K_0(O^{\Delta}) = K_0(O) = K_0(D^b(O)).$$

To $\Lambda^{\mathbf{d}}$ we associate simply the Serre subcategory $O^{\mathbf{d}}$ of $O$ generated by all $L(\lambda)$ with $\lambda \in \Lambda^{\mathbf{d}}$. More concretely, this is just a direct summand, specified by $k$, of the full subcategory $O^{\mathfrak{p}_{\mathbf{d}}}$ of $O$ of all modules which are locally finite over the standard parabolic subalgebra $\mathfrak{p}_{\mathbf{d}}$ with Levi factor $\mathfrak{gl}_{d_1} \oplus \cdots \oplus \mathfrak{gl}_{d_m}$. Sending $\vec{e} \in \bigwedge^{\mathbf{d}} V_{\pm}$ from (2.9) to the class of the $\mathfrak{p}_{\mathbf{d}}$-parabolic Verma module $\Delta^{\mathfrak{p}_{\mathbf{d}}}(\text{wt}(\vec{e}))$ (a standard object for the induced highest weight structure on $O^{\mathbf{d}}$) with highest weight $\text{wt}(\vec{e})$, defines an isomorphism of abelian groups

$$(2.11) \qquad (\bigwedge^{\mathbf{d}} V_{\pm}^{\mathbb{Z}}) \otimes_{\mathbb{Z}[v^{\pm 1}]} \mathbb{Z} \cong K_0(O^{\mathbf{d}}) = K_0(D^b(O^{\mathbf{d}})), \quad \vec{e} \mapsto [\Delta^{\mathfrak{p}_{\mathbf{d}}}(\text{wt}(\vec{e}))].$$

Here $\mathbb{Z}$ is a $\mathbb{Z}[v^{\pm 1}]$-module via $v \mapsto 1 \in \mathbb{Z}$ and $V_{\pm}^{\mathbb{Z}}$ denotes the $\mathbb{Z}[v^{\pm 1}]$-module in $V_{\pm}$ spanned by the $e_i$. We like to find functors realising the $U_v(\mathfrak{gl}_k)$-action and also incorporate $v$.

Tensoring with finite-dimensional representations of $\mathfrak{g}$ provides exact endofunctors of $O$. These functors and their direct summands form the monoidal category $\mathcal{P}_N$ of *projective functors*. Describing their effect on Verma modules is easy (but hard on other objects):

**Example 2.21.** If $U$ is a finite-dimensional representation of $\mathfrak{g}$, then $\Delta(\lambda) \otimes U \in O^{\Delta}$. The subquotients in a $\Delta$-flag are the $\Delta(\lambda + \nu)$, where $\nu$ runs through the multiset $P(U)$ of weights $\nu$ of $U$ with multiplicity $\dim U_{\nu}$. Thus, $[\Delta(\lambda) \otimes U] = \sum_{\nu \in P(U)} [\Delta(\lambda + \nu)]$ in $K_0(O^{\Delta})$. Important examples are $U = \mathbb{C}^N$ or $U = (\mathbb{C}^N)^*$, where $[\Delta(\lambda) \otimes U] = \sum_{i=1}^N [\Delta(\lambda \pm \delta_i)]$ respectively.

**Lemma 2.22.** *The functors* $E, F : O \mapsto O, M \mapsto M \otimes U$ *with* $U = \mathbb{C}^N$ *and* $U = (\mathbb{C}^N)^*$ *respectively, decompose into direct summands* $E = \bigoplus_{i \in \mathbb{Z}} E_i$, $F = \bigoplus_{i \in \mathbb{Z}} F_i$ *with*

$$(2.12) \quad [E_i \Delta(\lambda)] = \sum_{\{j \mid \lambda_j = i\}} [\Delta(\lambda + \delta_j)], \quad [F_i \Delta(\lambda)] = \sum_{\{j \mid \lambda_j = i+1\}} [\Delta(\lambda - \delta_j)].$$

This is an easy consequence of Example 2.21 and the fact that $O$ decomposes into summands $O_{\mathbf{c}}$ labelled by $S_N$-orbits $\mathbf{c}$ in $\mathbb{Z}^N$. Here $O_{\mathbf{c}}$ denotes the Serre subcategory of $O$ generated by $L(\lambda)$ with $\lambda \in \mathbf{c}$. Under (2.11), $O^{\mathbf{d}} \cap O_{\mathbf{c}}$ corresponds to a weight space. By definition, the functors $E_i$ and $F_i$ preserve $O^{\mathfrak{p}_{\mathbf{d}}}$, even $O^{\mathbf{d}}$ if $1 \leq i \leq k-1$. Formulas (2.12) resemble Lie algebra actions. Generalised formulas from Example 2.21 for $O^{\mathbf{d}}$ imply that the induced action on $K_0(O^{\mathbf{d}})$ agrees via (2.11) with the $(v \mapsto 1)$ specialised $U_v(\mathfrak{gl}_k)$-action on $\bigwedge^{\mathbf{d}} V_+^{\mathbb{Z}}$. (Note the positive sign $+$ here!).

**Remark 2.23.** Inside $\mathcal{P}_N$, the Hecke category appears naturally: It is known that $\mathcal{P}_N$ is the Karoubian closure of the additive monoidal category generated by $E_i$, $F_i$. The proof relies on a monoidal equivalence $\mathcal{P}_N \simeq \mathcal{HC}_N$ with a certain category $\mathcal{HC}_N$ of Harish–Chandra bimodules. Via Soergel's functor $\mathbb{V}$ from [132] and its extension in [136], $\mathcal{P}_N$ is equivalent to a category of *singular Soergel bimodules*. By restriction to endofunctors of $O_0 := O_{\mathbf{c}}$ with $0 \in \mathbf{c}$ one gets $f(\mathcal{SBim}_N) \simeq f(\mathcal{H}_N^{\mathrm{geo}})$ as a full monoidal subcategory, where $f$ means that we forget the grading. Remarkably, a classification of indecomposable projective functors for the categories $O^{\mathbf{d}}$ was only recently obtained [84], based on advances in the 2-representation theory of Hecke algebras, i.e. the representation theory of categorified Hecke algebras.

To incorporate $v$, we work with a *graded version* $\hat{O}$ of $O$ and its Serre subcategories as defined in [12], i.e. with graded modules over the endomorphism ring $A$ of a minimal projective generator of $O$ equipped with the Koszul grading from [12].

**Remark 2.24.** The origin of the grading is an equivalence of additive categories between $\mathrm{Proj}(O_0)$ and the full subcategory of $R_n$-mod of *Soergel modules* $\mathbb{C} \otimes_{R^w} M$ for $M \in f(\mathcal{SBim}_N)$ which has an obvious graded lift. We get a graded version of $\mathrm{Proj}(O_0)$ and then also of $O_0$. Note that $\mathcal{SBim}_N$ obviously acts on this category by tensoring over $R$ from the right. With some extra work, all Lie theoretic categories and functors used here can be lifted to a graded version. A general approach to lift modules (e.g. (parabolic) Verma modules) and the above functors to the graded setting is developed in [137].

**Lemma 2.25.** *Any choice of graded lift* $\hat{\Delta}^{\mathfrak{p}_{\mathbf{d}}}(\mathrm{wt}(\vec{e}))$ *of* $\Delta^{\mathfrak{p}_{\mathbf{d}}}(\mathrm{wt}(\vec{e}))$ *lifts* (2.11) *to an isomorphism of* $\mathbb{Z}[v^{\pm 1}]$-*modules* ($V_{\pm}^{\mathbb{Z}}$ *denotes the* $\mathbb{Z}[v^{\pm}]$-*submodule of* $V_{\pm}$ *spanned by the* $e_i$):

$$(2.13) \quad \Psi : \bigwedge^{\mathbf{d}} V_{\pm}^{\mathbb{Z}} \cong K_0(\hat{O}^{\mathbf{d}}) = K_0(D^b(\hat{O}^{\mathbf{d}})), \quad \vec{e} \mapsto [\hat{\Delta}^{\mathfrak{p}_{\mathbf{d}}}(\mathrm{wt}(\vec{e}))].$$

We realised now $\bigwedge^{\mathbf{d}} V_{\pm}^{\mathbb{Z}}$ as the Grothendieck group of a category and want to lift morphisms and the $U_v(\mathfrak{gl}_m)$-action (1.19) to functors. We first consider $\bigwedge^{\mathbf{d}} V_+^{\mathbb{Z}}$. If $\mathfrak{p}_{\mathbf{d}'} \subset \mathfrak{p}_{\mathbf{d}}$ are two standard parabolic subalgebras in $\mathfrak{gl}_N$, then there is the exact inclusion functor incl and its left adjoint *Zuckerman functor* $Z$ of taking the largest quotient in the target category, $\mathrm{incl} : O^{\mathbf{d}} \rightleftarrows O^{\mathbf{d}'} : Z$. Now incl and the derived functor $\mathcal{L}Z$ induce morphisms on $K_0$ which we connect to (1.14). Recall Proposition 1.7 and observe that $\mathcal{W}^+$ is generated as category by *basic webs* which look like a generator from (1.9) with identities to the left and right. To each basic web $t$ we associate a functor $\text{MSS}^+(t)$, which is, up to an overall shift, the obvious graded lift $\hat{\text{incl}}$ or $\mathcal{L}\hat{Z}$ of the inclusion respectively the derived Zuckerman functor (with hopefully self-explanatory notation):

$$(2.14) \qquad \curlyvee_{\mathbf{d}}^{\mathbf{d}'} := \hat{\text{incl}}[-ab]\langle -ab \rangle : \quad D^b(\hat{\mathcal{O}}^{\mathbf{d}}) \rightleftarrows D^b(\hat{\mathcal{O}}^{\mathbf{d}'}) \quad : \mathcal{L}\hat{Z} =: \curlywedge_{\mathbf{d}'}^{\mathbf{d}}$$

To each composition $t$ of basic web diagrams assign the composition $\text{MSS}^+(t)$ of functors.

**Example 2.26.** Let $k = 2$ and consider the webs (1.9) for $a = b = 1$ with induced morphisms $\wedge_{(2,0)}^{(1,1)} : \bigwedge^{(2,0)} V_+^{\mathbb{Z}} \rightleftarrows \bigwedge^{(1,1)} V_+^{\mathbb{Z}} : \wedge_{(1,1)}^{(2,0)}$. To $e_2 \wedge e_1$ we associate $\Delta^{p(2,0)}((2,1))$ which is just the trivial $\mathfrak{gl}_2$-module $\mathbb{C}$. The BGG resolution $\Delta^{p(1,1)}((1,2)) \to \boxed{\Delta^{p(1,1)}((2,1))}$ of $\mathbb{C}$ implies that $\text{incl}[-1]$ induces the map $e_2 \wedge e_1 \mapsto -v^{-1}(e_2 \otimes e_1 - v e_1 \otimes e_2)$ with $v = 1$ on the Grothendieck group. On the other hand, $\mathcal{L}Z\Delta^{p(1,1)}((2,1)) = Z\Delta^{p(1,1)}((2,1)) = \Delta^{p(2,0)}((2,1))$ and $\mathcal{L}Z\Delta^{(1,1)}((1,2)) = \Delta^{p(2,0)}((2,1))[-1]$ induce $e_2 \otimes e_1 \mapsto e_2 \wedge e_1, e_1 \otimes e_2 \mapsto -v$. Picking a lift $\hat{\Delta}^{p(2,0)}((2,1))$ with graded resolution $\hat{\Delta}^{p(1,1)}((1,2))\langle 1 \rangle \to \boxed{\hat{\Delta}^{p(1,1)}((2,1))}$ gives formulas (1.14).

The following summarises results from **[105]**, **[139]**, **[41]** and categorifies $q$-skew Howe duality:

**Theorem 2.27.** *Let $t$ be a basic web from $\mathbf{d}'$ to $\mathbf{d}$ with corresponding homomorphism $\Phi_+(t)$ from Proposition 1.7. Then there are choices of graded lifts in (2.13) and of (2.12), such that the following diagram commutes for $1 \leq i \leq k - 1$ (also for $E_i$ replaced by $F_i$),*

$$\begin{array}{ccccccccc}
& & & & \Psi & & & & \\
\bigwedge^{\mathbf{d}'} V_+^{\mathbb{Z}} & \xrightarrow{\Phi_+(t)} & \bigwedge^{\mathbf{d}} V_+^{\mathbb{Z}} & \xrightarrow{\Psi} & K_0(\hat{\mathcal{O}}^{\mathbf{d}}) & = & K_0(D^b(\hat{\mathcal{O}}^{\mathbf{d}})) \xleftarrow{[\text{MSS}^+(t)]} K_0(D^b(\hat{\mathcal{O}}^{\mathbf{d}'})) & = & K_0(\hat{\mathcal{O}}^{\mathbf{d}'}) \\
\downarrow E_i & & \downarrow E_i & & \downarrow [\hat{E}_i] & & \downarrow [\hat{E}_i] \qquad \downarrow [\hat{E}_i] & & \downarrow [\hat{E}_i] \\
\bigwedge^{\mathbf{d}'} V_+^{\mathbb{Z}} & \xrightarrow{\Phi_+(t)} & \bigwedge^{\mathbf{d}} V_+^{\mathbb{Z}} & \xrightarrow{\Psi} & K_0(\hat{\mathcal{O}}^{\mathbf{d}}) & = & K_0(D^b(\hat{\mathcal{O}}^{\mathbf{d}})) \xleftarrow{[\text{MSS}^+(t)]} K_0(D^b(\hat{\mathcal{O}}^{\mathbf{d}'})) & = & K_0(\hat{\mathcal{O}}^{\mathbf{d}'}) \\
& & & & \Psi & & & &
\end{array}$$

*Moreover, the family of functors $\hat{E}_i$, $\hat{F}_i$ naturally commutes with the functors $\text{MSS}^+(t)$ associated with webs. On the Grothendieck group they induce skew Howe duality (1.18)*

$$(2.15) \qquad U_v^{\mathbb{Z}}(\mathfrak{gl}_k) \quad \curvearrowright \quad X := \bigwedge^{\bullet}(V_+^{\mathbb{Z}}(k) \otimes V_-^{\mathbb{Z}}(m)) \quad \curvearrowleft \quad (U_v^{\mathbb{Z}}(\mathfrak{gl}_m))^{\text{op}}$$

*(with $\mathbb{Z}$ referring to Lusztig's integral version of the quantum group). The action of $D_j$ is hereby categorified by an appropriate grading shift on each categorified weight space.*

**Example 2.28.** We turned each summand from (1.19) into a category $D^{\mathbf{d}} := D^b(\hat{\mathcal{O}}^{\mathbf{d}})$, the $U_v(\mathfrak{gl}_k)$-action into functors $\hat{E}_i$, $\hat{F}_i$ and the action by $E$ and $-F$ into functors from (2.14):

$$\begin{array}{c}
D^{(2,1,0)} \qquad\qquad D^{(1,1,1)} \qquad\qquad D^{(0,1,2)} \\
{}^{\text{incl}[-2]\langle -2\rangle} \nearrow \qquad {}_{\mathcal{L}\hat{Z}=\text{id}} {}^{\text{incl}[-1]\langle -1\rangle} \nearrow \qquad {}_{\mathcal{L}\hat{Z}} \quad {}^{\text{incl}} \nearrow \qquad {}_{\mathcal{L}\hat{Z}} \\
D^{(3,0)} \quad \oplus \quad D^{(2,1)} \quad \oplus \quad D^{(1,2)} \quad \oplus \quad D^{(0,3)} \\
{}_{\mathcal{L}\hat{Z}} \searrow \quad {}^{\text{incl}} \quad {}_{\mathcal{L}\hat{Z}} \searrow \quad {}^{\text{incl}[-1]\langle -1\rangle} \quad {}_{\mathcal{L}\hat{Z}} \searrow \quad {}^{\text{incl}[-2]\langle -2\rangle} \\
D^{(2,1,0)} \qquad\qquad D^{(1,1,1)} \qquad\qquad D^{(0,1,2)}
\end{array}$$

In this categorified $q$-skew Howe duality, the two sides seem to be asymmetric. The action of $U_v^{\mathbb{Z}}(\mathfrak{gl}_k)$ is given by exact functors on the abelian categories, whereas $U_v^{\mathbb{Z}}(\mathfrak{gl}_m)$ acts by derived functors (note $V_+$ versus $V_-$). This asymmetry is explained nicely via Koszul (self-)duality [12], [106]. This directly gives an analogue of the theorem for $V_-$ instead of $V_+$.

**Remark 2.29.** *Koszul duality* means an equivalence $D^b(\hat{O}_{\mathbf{c}}^{\mathrm{pa}}) \simeq D^b(\hat{O}_{\mathbf{d}}^{\mathrm{pc}})$ which swaps the two types of functors [106]. Passing to the Grothendieck groups, it induces an isomorphism of groups $\bigwedge^{\bullet}(V_+^{\mathbb{Z}}(k) \otimes V_-^{\mathbb{Z}}(m)) \cong \bigwedge^{\bullet}(V_-^{\mathbb{Z}}(k) \otimes V_+^{\mathbb{Z}}(m))$. The parameters $v, \gamma, \eta$ from 1. III reflect important properties of this duality: it does not commute with grading shifts ($v \mapsto -v$ encoded by $\eta$) and does not preserve the standard $t$-structures [106] (encoded by $\gamma$).

Under Koszul duality, the derived functors $\mathrm{MSS}^+(t)$ from (2.14) turn into exact projective functors $\mathrm{MSS}^-(t)$ between the corresponding abelian categories. We use now these easier functors to construct tangle invariants with values in the homotopy categories $K^b(\hat{O}^{\mathbf{d}}, \hat{O}^{\mathbf{d}'})$ of exact functors from $\hat{O}^{\mathbf{d}}$ to $\hat{O}^{\mathbf{d}'}$. From [105] it follows that the relations in $\mathcal{W}^\eta$ can be interpreted in terms of isomorphisms of functors $\mathrm{MSS}^-(t)$. Thus we have cf. (1.8) an exact functor assigned to any basic tangle diagram except the crossings to which we assign the following complexes (possibly with identity strands added) given by canonical adjunction morphisms

$$(2.16) \quad \mathrm{MSS}^-\left(\diagup\!\!\!\diagdown\right): \left(\boxed{\mathrm{id}\langle 1\rangle} \to \mathrm{MSS}^-(\diagdown\!\!\diagup)\right)\langle -k\rangle, \quad \mathrm{MSS}^-\left(\diagup\!\!\!\diagdown\right): \left(\mathrm{MSS}^-(\diagdown\!\!\diagup) \to \boxed{\mathrm{id}\langle -1\rangle}\right)\langle k\rangle.$$

The following is proved in [105] and the Koszul dual version in [139]:

**Theorem 2.30.** *Let $t$ be an oriented tangle with a planar projection $t_1 \cdots t_r$ written in terms of basic tangle diagrams. Let $\mathrm{wRT}^-(t_1 \cdots t_r): \bigwedge^{\mathbf{d}} V_-^{\mathbb{Z}} \to \bigwedge^{\mathbf{d}'} V_-^{\mathbb{Z}}$. Then the composition*

$$(2.17) \quad \mathrm{MSS}^-(t_1) \cdots \mathrm{MSS}^-(t_r) \in K^b(\hat{O}^{\mathbf{d}}, \hat{O}^{\mathbf{d}'})$$

*is independent of the chosen projection. Thus, $t \mapsto \mathrm{MSS}^-(t)$ provides an invariant of oriented tangles. The induced morphism $K_0(D^b(O^{\mathbf{d}})) \to K_0(D^b(O^{\mathbf{d}'}))$ agrees via (2.13) with appropriate graded lifts with $\mathrm{wRT}^-(t)$. Analogously for $\mathrm{wRT}^+(t)$ using the Koszul dual functors.*

**Corollary 2.31.** *In case $t$ is a link, the categories $\hat{O}^{\mathbf{d}}, \hat{O}^{\mathbf{d}'}$ in (2.17) can be identified canonically with the category of graded vector spaces. Thus we obtain a bigraded link homology.*

**Remark 2.32.** Let $k = 2$. Then the invariant $\mathrm{MSS}^-$ was first defined in [135] based on [16], where it was observed that for non-quantised $\mathfrak{gl}_2$, the action of the Temperley-Lieb algebra on $(\mathbb{C}^2)^{\otimes n}$ can be categorified using category $O$. In [134] it is shown that the Khovanov complex for an oriented link agrees with the value of $\mathrm{MSS}^-$ by an explicit description of the involved categories as modules over (an extension of) *Khovanov's arc algebra*. Using [26] one can also match $\mathrm{MSS}^-$ with Khovanov's tangle invariant [79] via an equivalence of categories. Via [1], [99] which realises the extended arc algebra from [134] in terms of Fukaya-Seidel categories a rigorous categorical equivalence from $\mathrm{MSS}^-$ to the symplectic Khovanov invariant holds. A weaker combinatorial identification, the equality of the bigraded homology groups, is established (in fact for all known algebraic-geometric link homologies) in [97].

**Remark 2.33.** For $k = 2$, functoriality (as in (2.1)) of $\text{MSS}^\pm$ is reduced to that of $\text{Kh}_{\text{sgn}}$. For general $k$ we expect functoriality to follow from the functoriality results in **[44]**.

**Remark 2.34.** We focused here on defining the involved functors and describing their action on the Grothendieck group, although all defining relations in the quantum group or web category can in fact be turned into actual relations, i.e. isomorphisms, between functors.

**Remark 2.35.** Formula (1.15) is implicitly categorified via $\text{MSS}^-(\beta_{a,b})$: For $a = b = 1$ this holds by (2.16). A composition of those, cf. illustration in (3.1), gives the braiding morphism for $W^{\otimes a} \otimes W^{\otimes b}$ and restriction to $\hat{\mathcal{O}}^{(a,b)}$ then $\text{MSS}^-(\beta_{a,b})$. One can verify purely combinatorially based on **[105]**, that the complex $\text{MSS}^-(\beta_{a,b})$ of exact functors can be written with entries encoded by (1.15). Lie theoretically $\text{MSS}^-(\beta_{a,b})$ is easy to describe as the derived functor of a classical shuffling functor **[135]**, **[105]** which gets *reinterpreted* in terms of the explicit complexes. This is opposite to most categorifications, where the braiding is *defined* by explicit complexes indicated by (1.15), e.g. **[145]**, **[115]**. The construction of categorified braid group actions from categorified Lie algebra actions using Rickard complexes goes back to **[32]**.

**Remark 2.36.** Categorifications of (parts of) $q$-skew Howe duality were obtained and used in many ways in recent years. The significance of the above construction is the fact that *both* quantum group actions are visible. This is in particular not the case in diagrammatic or foam based categorifications, since a (diagrammatic) replacement of the derived functors is missing. It would be nice to find a general theory towards categorifications of dualities, in particular for those in Remark 1.12 where a categorification so far only exists for iv) **[41]**.

**Q4**: *Are there other interesting Koszul self-dual categories? What do they categorify?*

### 2.3.1. Towards 2-Representation Theory: Categorical actions.

Two basic questions arise from the above construction: is there a conceptual source for isomorphisms specifying the desired relations between the functors (topologically speaking the values for tangle cobordisms)? To which extent are such categorifications unique? Both questions are addressed with the concept of *categorical Lie algebra actions* **[32]**, **[125]**, which we try to motivate based on our example. The categorified quantum groups due to Khovanov–Lauda **[81]** and Rouquier **[125]** occur in this context naturally. Adjunction morphisms between functors are used to specify commutation relations (e.g. for $\text{E}_i$ and $\text{F}_i$) and most tangle cobordisms. More involved are the Serre relations between the $\text{E}_i$ which arise from endomorphisms (= natural transformations) of (powers of) $\text{E} = \_ \otimes \mathbb{C}^N$ which we construct below. In the general construction **[81]**, **[125]** such morphisms constitute 2-morphisms of a 2-category.

An obvious choice for $\mathfrak{s} \in \text{End}(\text{E}^2)$ is the *flip* morphism given on $\mathfrak{g}$-modules $M$ by $\mathfrak{s}_M : M \otimes \mathbb{C}^N \otimes \mathbb{C}^N \to M \otimes \mathbb{C}^N \otimes \mathbb{C}^N$, $m \otimes u_1 \otimes u_2 \mapsto m \otimes u_2 \otimes u_1$. The action maps $\mathbb{x}_M : M \otimes \mathbb{C}^N \to M \otimes \mathbb{C}^N$ of the *Casimir element* $\Omega = \sum_{i,j=1}^N \text{E}_{i,j} \otimes \text{E}_{j,i} \in \mathfrak{g} \otimes \mathfrak{g}$ define an endomorphism $\mathbb{x} \in \text{End}(\text{E})$. More generally define endomorphisms $\mathbb{x}_j, \mathfrak{s}_j$ of $\text{E}^n$, $n \geq 0$, via

(2.18) $\quad (\mathbb{x}_j)_M := \text{E}^{n-j}(\mathbb{x}_{\text{E}^{j-1}(M)}) \quad \text{and} \quad (\mathfrak{s}_j)_M := \text{E}^{n-j-1}(\mathfrak{s}_{\text{E}^{j-1}(M)}).$

One easily verifies that these endomorphisms satisfy the defining relations of a *degenerate affine Hecke algebra* $\text{H}_n^{\text{daff}}$. That means the $\mathbb{x}_i$ commute (defining a subalgebra

23        Categorification: tangle invariants and TQFTs

$\mathbb{C}[\mathbbm{x}_1, \ldots, \mathbbm{x}_n])$, the $\mathbbm{s}_i$ satisfy the Coxeter relations of the symmetric groups (defining a subalgebra $\mathbb{C}[S_n]$) and the two sets of generators interact via the degenerate semidirect product relations $\mathbbm{s}_j \mathbbm{x}_{j+1} = \mathbbm{x}_j \mathbbm{s}_j + 1$ and $\mathbbm{s}_j \mathbbm{x}_l = \mathbbm{x}_l \mathbbm{s}_j$ for $l \neq j, j+1$.

Amazingly (although easy to verify with Example 2.21), $E_i$ from (2.12) equals the *(generalised) $i$-eigenspace subfunctor* for $\mathbbm{x}$ of E, i.e. $E_i(M) = \sum_{l \geq 0} \ker((\mathbbm{x}_M - i)^l)$. To (re)define $F_i$ consistently we use that F is *right* adjoint to E. With a *fixed* counit $c : EF \to \text{id}$ and unit $c^* : \text{id} \to FE$ we define elements $\mathbbm{x}' \in \text{End}(F)$ and $\mathbbm{s}' \in \text{End}(F^2)$ following [125]:

$$(2.19) \qquad \mathbbm{x}' := F(c) \circ F(\mathbbm{x})_F \circ c^*_F, \quad \mathbbm{s}' := F^2(c) \circ F^2 E(c)_F \circ F^2(\mathbbm{s})_{F^2} \circ E(c^*)_{EF^2} \circ c_{F^2}.$$

Then F inherits a decomposition into $i$-eigenspace functors $F_i$, $1 \leq i \leq k$, for $\mathbbm{x}'$. By [125] the *bi*adjointness of $E_i$ and $F_i$ follows. Thus, these functors are exact, send projectives to projectives and provide a *based categorification*, i.e. induce on the Grothendieck group $K_0^\oplus(\text{Proj}(\mathcal{O}^\mathbf{d})) \otimes_\mathbb{Z} \mathbb{C}$ the structure of an integrable $\mathfrak{sl}_k$-module $\bigwedge^\mathbf{d} \overline{V}_+$, where the classes of the indecomposable projectives are (a basis of) weight vectors. To agree with the existing literature, we work now with $\mathfrak{sl}_k$. These ingredients and properties were axiomatised in [125]:

**Definition 2.37.** Let $\mathfrak{C}$ be a $\mathbb{C}$-linear abelian finite length category with enough projectives. Then, a *categorical* $\mathfrak{sl}_k$-*action* on $\mathfrak{C}$ (categorifying $C$) consists of

- an endofunctor E with a right adjoint F specified by a counit $c$ and a unit $c^*$
- an element $\mathbbm{s} \in \text{End}(E^2)$ and an endomorphism $\mathbbm{x} \in \text{End}(E)$

which satisfy: $E = \oplus_{i=1}^k E_i$, where $E_i$ is the $i$-eigenspace subfunctor for $\mathbbm{x}$, the endomorphisms $\mathbbm{s}_j$, $\mathbbm{x}_j$ defined via (2.18) satisfy the relations of $H_n^{\text{daff}}$, the functor F is right adjoint to E, and finally with the definition of $F_i$ as $i$-eigenspace functor for $\mathbbm{x}'$ as in (2.19), the functors $E_i$ and $F_i$ define a based categorification of an integrable $\mathfrak{sl}_k$-module $C$.

We get categorifications of the tensor products $\bigwedge^\mathbf{d} \overline{V}_+$ of $\mathfrak{sl}_k$ exterior powers:

**Theorem 2.38.** *The constructions* (2.18), (2.11) *define a categorical* $\mathfrak{sl}_k$-*action on* $\mathcal{O}^\mathbf{d}$.

**Remark 2.39.** Definition 2.37 is the easiest example from the theory of categorical actions of Kac–Moody Lie algebras $\mathfrak{g}_{\text{KM}}$ [125], [81]. The degenerate affine Hecke algebra gets replaced by a more general *quiver Hecke algebra* (or KLR algebra) which is used to define a certain graded 2-category $^2\dot{U}_v(\mathfrak{g}_{\text{KM}})$ categorifying $\dot{U}_v(\mathfrak{g}_{\text{KM}})$ cf. Remark 1.11. The definition is via generators and relations, algebraically [125] or diagrammatically [81], and matched in [21].

**Application 2.40.** A nice situation occurs when the morphisms (2.18) generate the endomorphism ring of an object $E^n(M)$ and the kernel is controlled by a cyclotomic quotient of $H_n^{\text{daff}}$. Then this ring can be determined explicitly. If moreover every indecomposable projective object in $\mathfrak{C}$ arises as a summand of $E^n(M)$ for $n \gg 0$, one might construct equivalences by determining and matching endomorphism rings of projective generators instead of providing a functor. This idea is applied e.g. in [25] to the category $\mathcal{F}(a|b)$ of finite-dimensional representations of the linear supergroup $\text{GL}(a|b)$: $\mathcal{F}(a|b)$ is equivalent to the category of modules for an infinite-dimensional analogue of Khovanov's arc algebra, Remarks 2.2, 2.32. The notion *higher Schur–Weyl duality* [22] formalises such nice Lie theoretic situations.

Definition 2.37 significantly rigidifies the involved category $\mathfrak{C}$. If $C$ is finite dimensional irreducible of highest weight $\xi$, its weight space decomposition implies a decomposition of $\mathfrak{C}$ into direct summands $\mathfrak{C}_\lambda$ (cf. with $O^{\mathbf{d}}$) and the *Uniqueness Theorem*, a very special case of Rouquier's *Universality Theorem*, holds [32], [125]: a *minimal* (i.e. $\mathfrak{C}_\xi \simeq$ Vect) categorification of such $C$ is unique up to *strong equivalence* meaning an equivalence of categories $\Gamma : \mathfrak{C} \to \mathfrak{C}'$ with an isomorphism $\phi : \Gamma \mathrm{E} \cong \mathrm{E} \Gamma$ satisfying the expected compatibilities with $\mathrm{x}$, $\mathrm{s}$. Uniqueness allows to establish abstractly equivalences of categories.

**Application 2.41.** The Uniqueness Theorem is powerful even for $\mathfrak{sl}_2$-modules. It is used in [32] to prove Broué's abelian defect group conjecture for symmetric groups, one of the most famous conjectures in modular representation theory of finite groups.

**Application 2.42.** By the Universality Theorem, the $k \geq 3$ generalisations [96], [98] of Khovanov's arc algebras are Morita equivalent to certain cyclotomic quotients of $\mathrm{H}_n^{\mathrm{daff}}$. These algebras should provide an algebraic construction of the $\mathtt{MSS}^-$ invariants as in Remark 2.32.

**2.3.2. Tensor product categorifications.** The Uniqueness Theorem heavily relies on the fact that finite-dimensional irreducible modules are generated by their highest weight vectors and thus does not directly apply to Theorem 2.38. A general theory for *the process of taking* tensor products of categorifications is still missing. The *naive* outer tensor product of categorical $\mathfrak{sl}_k$-actions has the desired $\mathrm{K}_0$, but only a categorical $\mathfrak{sl}_k \oplus \mathfrak{sl}_k$-action. Losev and Webster [92] give an axiomatic definition of categorification of a *given* tensor product.

Their definition uses the *reverse dominance ordering* on weights in a tensor product. Concretely consider again the $\mathfrak{sl}_k$-module $\bigwedge^{\mathbf{d}} \overline{V}_+$ with $\mathbf{d} = (d_1, \ldots, d_m)$. View its weights as tuples $\lambda = (\lambda_1, \ldots, \lambda_m)$ of $\mathfrak{sl}_k$-weights, and $\lambda \geq \mu$ if $\lambda_1 + \cdots + \lambda_m = \mu_1 + \cdots + \mu_m$ and $\lambda_1 + \cdots + \lambda_i \leq \mu_1 + \cdots + \mu_i$ (in the usual ordering on $\mathfrak{sl}_k$-weights) for each $i < m$. Via (2.9), this ordering translates into the Bruhat ordering on $\Lambda^{\mathbf{d}}$. The following is a reformulation of the original definition [92] following closely [23].

**Definition 2.43.** A *tensor product categorification* of the $\mathfrak{sl}_k$-module $\bigwedge^{\mathbf{d}} \overline{V}_+$ is the same data as in Definition 2.37, but with the last property on based categorification replaced by

- $\mathfrak{C}$ is a highest weight category with respect to the poset $\Lambda^{\mathbf{d}}$
- the exact functors $\mathrm{E}_i$ and $\mathrm{F}_i$ send objects with $\Delta$-flags to objects with $\Delta$-flags
- under the isomorphism $\bigwedge^{\mathbf{d}} \overline{V}_+ \cong \mathrm{K}_0(\mathfrak{C}) \otimes_{\mathbb{Z}} \mathbb{C}$, $\vec{e} \mapsto [\Delta^{\mathfrak{C}}(\mathrm{wt}(\vec{e}))]$ as in (2.11), the actions of $E_i$ and $F_i$ correspond to the actions of $[\mathrm{E}_i]$ and $[\mathrm{F}_i]$, respectively.

**Theorem 2.44.** *The highest weight category $O^{\mathbf{d}}$ defines with the data from Theorem 2.38, the poset $\Lambda^{\mathbf{d}}$, (2.12) and (2.11) a tensor product categorification of the $\mathfrak{sl}_k$-module $\bigwedge^{\mathbf{d}} \overline{V}_+$.*

By the following result from [92] this is the only one up to strong equivalence:

**Theorem 2.45.** *A tensor product categorification of the $\mathfrak{sl}_k$-module $\bigwedge^{\mathbf{d}} \overline{V}_+$ is unique.*

**Remark 2.46.** Definition 2.43 is again a special case of a more general definition [92] which works for any Kac–Moody Lie algebra $\mathfrak{g}_{\mathrm{KM}}$ instead of $\mathfrak{sl}_k$ and any integrable highest weight module of $\mathfrak{g}_{KM}$ for each tensor factor. It requires the following adjustments. On the one hand, the action of the degenerate affine Hecke algebra gets replaced by a quiver Hecke algebra from Remark 2.39 or an even more general *Webster algebra* [142]. On the other hand, the highest weight category gets replaced by a *fully stratified* category. For the general theory of such generalisations of highest weight categories see [24].

**Remark 2.47.** Let us return to the graded setting to obtain categorifications of the $U_v(\mathfrak{gl}_k)$-modules $\bigwedge^{\mathbf{d}} V_+$. One can turn $\hat{O}^{\mathbf{d}}$ into a graded additive $\mathbb{C}$-linear 2-category. For this note that $\mathrm{Hom}_{O^{\mathbf{d}}}(M, N) = \bigoplus_{j \in \mathbb{Z}} \mathrm{Hom}_{\hat{O}^{\mathbf{d}}}(\hat{M}, \hat{N}\langle j \rangle)$ for $M, N$ which have graded lifts $\hat{M}, \hat{N}$, similarly for functors. Objects in this 2-category are weights $\mathbf{c}$ of $\bigwedge^{\mathbf{d}} V_+$, but thought of as the corresponding summands in $\hat{O}^{\mathbf{d}}$ via (2.13). Morphisms are generated by *i*) the functors $\mathbb{1}_{\mathbf{c}}$ which are the identity on $\mathbf{c}$ and zero otherwise, *ii*) the functors $\hat{\mathsf{E}}_i \mathbb{1}_{\mathbf{c}}$ viewed as morphisms from $\mathbf{c}$ to $\mathbf{c} + \alpha_i$ where $\alpha_i$ is the corresponding simple root for $\mathfrak{gl}_k$ and *iii*) fixed right adjoints of *ii*) which are the $\mathbb{1}_{\mathbf{c}} \hat{\mathsf{F}}_i$ up to shifts. The 2-morphisms are generated by the homogeneous components of the natural transformations (2.18). From [125] it follows that this data defines a *(strong) 2-representation* of $\mathfrak{sl}_k$. By [29] it extends to an *action of* $^2\dot{U}_v(\mathfrak{q})$, called a *2-representation of* $\dot{U}_v(\mathfrak{q})$, for $\mathfrak{q} = \mathfrak{sl}_k$ and also for $\mathfrak{gl}_k$ by adding a grading shifting operator.

**Remark 2.48.** The original definition in [92] connects Theorem 2.44 with naive outer tensor products: for $\lambda \in \Lambda^d$ there are Serre subcategories $O^{\mathbf{d}}[<\lambda] \subset O^{\mathbf{d}}[\leq \lambda]$ in $O^{\mathbf{d}}$ generated by all $L(\mu)$ with $\mu < \lambda$ respectively $\mu \leq \lambda$. The *associated graded* $\bigoplus_{\lambda} O^{\mathbf{d}}[\leq \lambda] / O^{\mathbf{d}}[<\lambda]$ of $O^{\mathbf{d}}$ formed from the subquotients can be identified with the naive tensor product of the categorifications of the factors $\bigwedge^{d_i} \overline{V}_+$, see [127] for an explicit identification. The highest weight structure, explicitly the poset $\Lambda^{\mathbf{d}}$, creates a desired asymmetry, namely the asymmetry in the tensor factors (1.4) when passing to the graded/quantised setting as in Remark 2.47.

**Application 2.49.** Categorical actions are often used in (specifically modular and super) representation theory to create interesting gradings or to determine decomposition numbers. We sketch an example directly connected to our framework. In [23], tensor product categorifications were defined for the limit Lie algebra $\mathfrak{sl}_{\mathbb{Z}}$ and constructed for $M := \overline{V}_{\infty}^{\otimes a} \otimes (\overline{V}_{\infty}^*)^{\otimes b}$, very similar to above, using category $O$ for the Lie superalgebra $\mathfrak{gl}_{a|b}$. Here, $\overline{V}_{\infty} = \mathbb{C}^{\mathbb{Z}}$ is the natural representation of $\mathfrak{sl}_{\mathbb{Z}}$ and $\overline{V}_{\infty}^*$ its restricted dual. The basis vectors in $\mathbb{C}^{\mathbb{Z}}$ labelled by a length $k$ interval in $\mathbb{Z}$ span an $\mathfrak{sl}_k$-module $\overline{V}_+^{\otimes a} \otimes (\bigwedge^{k-1} \overline{V}_+)^{\otimes b}$. Theorems 2.44 and 2.45 allow to translate properties from $O(\mathfrak{gl}_{a+(k-1)b})$ to the super side [23]. This implies that the (integral blocks of) category $O$ for $\mathfrak{gl}(a|b)$ and the category $\mathcal{F}(a|b)$ of finite-dimensional representation of $GL(a|b)$ can be equipped with a Koszul grading. Moreover, the graded decomposition numbers are given by parabolic Kazhdan–Lusztig polynomials. In case of $\mathcal{F}(a|b)$ this grading agrees with the explicit construction in [25] from Application 2.40. For a generalisation to the more complicated orthosymplectic supergroups see [42].

### 2.4. Ad III: Categorified coloured tangle invariants and projectors

The coloured framed tangle invariant RT from Remark 1.4 involves ultimately tensor products of arbitrary finite dimensional irreducible $U_v(\mathfrak{gl}_k)$-modules, not only exterior pow-ers. A categorification of all such tensor products exists for $\mathfrak{sl}_k$ and $U_v(\mathfrak{sl}_k)$ [53], [142], [127], and by [142] even for any simple complex Lie algebra $\mathfrak{g}_s$. Webster [142] also ensures the existence of tensor product categorifications for $\mathfrak{g}_s$. He uses categories of graded modules over graded algebras which generalise quiver Hecke and quiver Schur algebras. These algebras are defined diagrammatically, so that all calculations are elementary, but usually not easy. The grading allows to get categorifications of $U_v(\mathfrak{g}_s)$-modules as in Remark 2.47 [142] with a direct generalisation of Theorem 2.30 to arbitrary $\mathfrak{g}_s$. Instead of formulating this in detail we indicate phenomena which occur, even for $\mathfrak{g}_s = \mathfrak{sl}_k$, when passing from tensor products of fundamental representations to arbitrary irreducible ones. On the way we construct tensor product categorifications for $\mathfrak{sl}_k$ using the results from the previous section. The $\mathfrak{sl}_k$-action extends by construction to a $\mathfrak{gl}_k$-action, even to a 2-representation of $\dot{U}_v(\mathfrak{gl}_k)$ when invoking gradings. However not all irreducible $\mathfrak{gl}_k$-modules occur in this way as tensor factors.

Any irreducible finite-dimensional $\mathfrak{sl}_k$-module is a quotient of some $\bigwedge^{\mathbf{d}} \overline{V}_+$ as in Section 2.3.2 such that its highest weight is the sum of the highest weights of the tensor factors. Taking tensor products $\bigwedge^{\mathbf{d}_{(1)}} \overline{V}_+ \otimes \cdots \otimes \bigwedge^{\mathbf{d}_{(r)}} \overline{V}_+ \twoheadrightarrow \overline{V}(\xi_1) \otimes \cdots \otimes \overline{V}(\xi_r) =: \overline{V}(\xi)$ realises finite tensor products $\overline{V}(\xi)$ also as quotients of some $\bigwedge^{\mathbf{d}} \overline{V}_+$ which we consider now. Via (2.11) the irreducible objects in $\mathcal{O}^{\mathbf{d}}$ give rise to a *special basis* (in fact the specialised Lusztig dual canonical basis) of $\bigwedge^{\mathbf{d}} \overline{V}_+$. It turns out that the kernel of the quotient to $\overline{V}(\xi)$ is spanned by a subset of these special basis vectors. Fix $1 \leq j \leq r$. Combinatorially, one can label standard basis vectors (1.12) in $\bigwedge^{\mathbf{d}_{(j)}} \overline{V}_+$ canonically by column strict tableaux and then basis vectors in $\overline{V}(\xi_j)$ by the set $I_j$ of semistandard tableaux, i.e. those which are additionally weakly row strict. The shape is determined by $\mathbf{d}_{(j)}$ or equivalently $\xi_j$ and fillings are from $\{1, \ldots, k\}$. Consider now the standard basis vectors $\vec{e}$ as in (2.9) which correspond to $m$-tuples not in $I_1 \times \cdots \times I_r$. They define (by taking the irreducible quotient of the corresponding parabolic Verma module in (2.13)) a set of irreducible objects $L(\text{wt}(\vec{e})) \in \mathcal{O}^{\mathbf{d}}$, thus a Serre subcategory $\mathcal{S}_\xi$ in $\mathcal{O}^{\mathbf{d}}$. We obtain a categorification of $\overline{V}(\xi)$ [127], implicitly [142]:

**Theorem 2.50.** *The Serre quotient $\mathcal{O}^{\mathbf{d}}/\mathcal{S}_\xi$ inherits a categorical $\mathfrak{sl}_k$-action from $\mathcal{O}^{\mathbf{d}}$. This is a tensor product categorification in the sense of [92] categorifying $\overline{V}(\xi)$ with the ordering on the labelling set of irreducible objects induced from $\Lambda^{\mathbf{d}}$.*
*From $\hat{\mathcal{O}}^{\mathbf{d}}$ as in Remark 2.47, the graded version $\hat{\mathcal{O}}^{\mathbf{d}}/\hat{\mathcal{S}}_\xi$ inherits an action of $^2\dot{U}_v(\mathfrak{gl}_k)$.*

The quotient functors $\pi_\xi : \hat{\mathcal{O}}^{\mathbf{d}} \to \hat{\mathcal{O}}^{\mathbf{d}}/\hat{\mathcal{S}}_\xi$ are exact and induce $\mathbb{Z}[v^{\pm 1}]$-linear morphism on the Grothendieck groups (which however usually do not split over $\mathbb{Z}[v^{\pm 1}]$):

(2.20) $$\begin{array}{ccc} \hat{\mathcal{O}}^{\mathbf{d}} & \xrightarrow{\pi_\xi} & \hat{\mathcal{O}}^{\mathbf{d}}/\hat{\mathcal{S}}_\xi \\ \downarrow{\scriptstyle K_0} & & \downarrow{\scriptstyle K_0} \\ \bigwedge^{\mathbf{d}} V_+^{\mathbb{Z}} & \xrightarrow{[\pi_\xi]} & V(\xi)^{\mathbb{Z}} \end{array}$$

(2.21) $(V_{+,\mathfrak{gl}_2}^{\mathbb{Z}})^{\otimes m} \xrightarrow[\text{split ?}]{[\pi_m]:=[\pi_\xi]} V_+^{\mathbb{Z}}(\xi), \quad \xi := m\omega_1$

where $\omega_1$ is the first fundamental weight
*Jones-Wenzl projector:* $\text{JW}_m := [\mathcal{L}\iota_m] \circ [\pi_m]$

**Remark 2.51.** Theorem 2.50 requires a more general version of Definition 2.43 from [92], see Remark 2.46, as the quotient category $\mathcal{O}^{\mathbf{d}}/\mathcal{S}_\xi$ might not be highest weight, but only fully stratified. Combinatorially this is reflected in higher dimensional weight spaces of the tensor factors of $\overline{V}(\xi)$. Using only fundamental representations avoids this problem as they are minuscule, and also avoids taking duals or inverses of determinant representations.

In contrast to $\hat{O}^{\mathbf{d}}$, the quotients $O^{\mathbf{d}}/\mathcal{S}_\xi$ usually have *infinite global dimension*. Thus, the computation of a derived left adjoint $\mathcal{L}\iota_\xi$ to the quotient functor $\pi_\xi$ requires *infinite* resolutions and *unbounded* derived categories. This becomes relevant in categorifications of coloured tangle invariants following the knot theoretic colouring via *cabling* and *projectors*. The idempotent functor $\mathrm{pr}_\xi := \mathcal{L}\iota_\xi \pi_\xi$ is a categorified projector.

**Remark 2.52.** Working with infinite complexes is delicate, in particular when Grothendieck groups or Euler characteristics are involved. To avoid an Eilenberg swindle and the collapse of the Grothendieck groups we work in the graded setting with certain subcategories $D^\nabla(\hat{O})$ of the unbounded derived category, such that $K_0(D^\nabla(\hat{O})) \cong K_0(\hat{O}) \otimes_{\mathbb{Z}[v^{\pm 1}]} \mathbb{Z}((v))$, see [2] for a precise definition. The functors $\pi_\xi$, $\mathcal{L}\iota_\xi$ induces then $\mathbb{Z}((v))$-linear maps

$$(2.22) \qquad [\pi_\xi] : (V_+^{\mathbb{Z}})^{\otimes m} \otimes_{\mathbb{Z}[v^{\pm 1}]} \mathbb{Z}((v)) \rightleftarrows V(\xi)^{\mathbb{Z}} \otimes_{\mathbb{Z}[v^{\pm 1}]} \mathbb{Z}((v)) : [\mathcal{L}\iota_\xi].$$

**Example 2.53.** In case of $U_v(\mathfrak{gl}_2)$, there is the quotient map $[\pi_2]$ as in (2.21) to the biggest irreducible quotient. Explicitly for $m = 2$, the following elements form a basis of the quotient

$$b_1 := [\pi_2](e_1 \otimes e_1), \quad b := [\mathrm{p}_2](e_1 \otimes e_2) = v^{-1}[\mathrm{p}_2](e_2 \otimes e_1), \quad b_2 := [\mathrm{p}_2](e_2 \otimes e_2).$$

A split of $[\mathrm{p}_2]$ over $\mathbb{C}(v)$ is given by $b_i \mapsto e_i \otimes e_i$ and $b \mapsto \frac{1}{[2]}(e_2 \otimes e_1 + v^{-1}e_1 \otimes e_2)$. Interpreting the latter as $(ve_2 \otimes e_1 + e_1 \otimes e_2)(1 - v^2 + v^4 - \cdots) \in (V_+^{\mathbb{Z}})^{\otimes 2} \otimes_{\mathbb{Z}[v^{\pm 1}]} \mathbb{Z}((v))$, we obtain the morphism $[\mathcal{L}\iota_2]$ induced via (2.22) from $\mathcal{L}\iota_2$ (without explicitly constructing $\mathcal{L}\iota_2$).

**Remark 2.54.** In case of $U_v(\mathfrak{gl}_2)$, the projector (2.21) from $V_{+,\mathfrak{gl}_2}^{\otimes m}$ onto the biggest irreducible summand is the *Jones–Wenzl projector* $\mathrm{JW}_m$. Categorified is this the first time in [53] using Serre quotient, and its (Koszul dual) $V_{-,\mathfrak{gl}_2}$-version in [34], [126] using Bar-Natan's approach to Khovanov homology respectively iterated categorified full twists.

As a special case, the RT value of the unknot coloured by the $U_v(\mathfrak{gl}_k)$-module $V^{\mathbb{Z}}(\xi)$ as in (2.20) can be categorified by taking the $\mathrm{MSS}^+$-value of $\sum_i d_i$ nested cups (viewed as a derived functor) followed by $\mathrm{pr}_\xi$ and followed by the value of $\sum_i d_i$ nested caps:

(2.23) the projector $\mathrm{pr}_\xi$ $\begin{array}{c}\mathcal{L}\iota_\xi \\ \pi_\xi\end{array}$ (2.24) the categorified value of the coloured unknot $\begin{array}{c}\cdots \quad \cdots \\ \mathcal{L}\iota_\xi \\ \pi_\xi \\ \cdots \quad \cdots\end{array}$

**Example 2.55.** In case $k = 2$ and $V(\xi)$ is 3-dimensional, the value of (2.24) can be realised as a complex in the unbounded homotopy category $K^-(\mathbb{C}\text{-}\mathrm{mod}^{\mathbb{Z}(v)})$. A lengthy calculation gives the graded Poincaré polynomial $v^2 t^2 + 1 + v^{-2} + \frac{v^{-6}t^{-2}(1+t^{-1})}{1-t^{-2}v^{-4}} \in \mathbb{Z}[t^{\pm 1}]((v))$ [138]. Its Euler characteristic equals $v^2 + 1 + v^{-2} = [3]$ which is indeed the $\mathrm{RT}_{V(\xi)}$-value of the unknot.

By a uniqueness result of categorified Jones–Wenzl projectors from [34], the value of the unknot from (2.24) or from Theorem 2.56 agrees with the Cooper–Krushkal categorified value [34] of the coloured unknot up to Koszul duality (i.e. a transformation $v \mapsto t^{-1}v^{-1}$).

Let $L$ be an oriented link with planar projection $D$ and colouring col assigning some $\overline{V}(\xi_c)$ to each components $c$ of $D$. Assume $(V_+^{\mathbb{Z}})^{\otimes m_c} \twoheadrightarrow V^{\mathbb{Z}}(\xi_c)$ as in (2.20). To $D$ we attach its *colour-cabled version* $D_{cc}$: we first replacing each strand in a component $c$ by its cabling, i.e. by $m_c$ parallel strands oriented as before. Then we write the result as a composition $t_1 \cdots t_r$ of basic tangle diagrams and finally place for one upwards pointing original strand in $D$ a projector (2.24) on its cabling. Let $\text{MSS}^+(D_{cc})$ be the associated composition of derived functors given by $\text{MSS}^+$ with additionally $\text{pr}_\xi$ included when the projector occurs. Using the identification from Corollary 2.31 we can apply this functor to the vector space $\mathbb{C}$ concentrated in bidegree zero to get an object $\text{MSS}^+(D_{cc})(\mathbb{C})$ in $D^\triangledown(\mathbb{C}\text{-mod}^{\mathbb{Z}(v)}) \subset K^-(\mathbb{C}\text{-mod}^{\mathbb{Z}(v)})$.

The categorification [138] of the coloured framed oriented tangle invariant with colours irreducible $U_v(\mathfrak{sl}_k)$-modules (or their $U_v(\mathfrak{gl}_k)$-versions (2.20)) implies for links:

**Theorem 2.56.** *The assignment $D \mapsto \text{MSS}^+(D_{cc})(\mathbb{C})$ defines an invariant $\text{MSS}^+_{col}$ of coloured framed oriented links. It induces on $K_0$ the coloured RT-invariant from Remark 1.4 for $\mathfrak{sl}_k$.*

**Remark 2.57.** The coloured knot invariants from Theorem 2.56 are usually infinite complexes, even for the unknot. The Poincaré series of $\text{MSS}^+(D_{cc})(\mathbb{C})$ has values in $\mathbb{Z}[t^{\pm 1}]((v))$. This is similar to the HHH-invariant, but we believe it is even harder to compute. For $k = 2$, these invariants should be directly connected to the invariants constructed in [31], where impressive explicit examples are computed. The occurring infinite series are secretly rewriting quotients $\frac{[a]}{[b]}$ of quantum numbers, see Example 2.53. A realisation of such quotients as Euler characteristic of an infinite complex is called *fractional Euler characteristic* in [54].

Recall from Section 1 that the HOMFLY-PT polynomial recovers the quantum $\mathfrak{gl}_k$ link invariants $\text{RT}_{V_-}$ by specialisation of **a** to $v^k$. One might expect a similar connection for the categorifications, i.e. between the triply graded KR link homology which by Theorem 2.12 is a categorification of the HOMFLY-PT-polynomial and $\text{MSS}^-$. Naive specialisation does not work, but there is a spectral sequence connecting the two theories, predicted in [40] and established in [117]. Also, recall from Section 2.2 that the approach to compute $\text{KR}(t_{(n,q)})$ and its limit $\text{KR}(t_{(n,\infty)})$ for torus links uses categorified projectors. On the other hand, $\text{MSS}^+(\bigcirc_{cc})(\mathbb{C})$, or its Koszul dual version $\text{MSS}^-(\bigcirc_{cc})(\mathbb{C})$, can be seen as a categorification of the closure of a projector. One again might expect a connection between $\text{KR}(\bigcirc)^n$ from Example 2.9, $\text{KR}(t_{(n,\infty)})$, and $\text{MSS}^\pm(\bigcirc_{cc})(\mathbb{C})$. The following reformulates conjectures from [60]:

**Conjecture 2.58.** *The algebra $B = \mathbb{C}[u_1, \ldots, u_n] \otimes \bigwedge^\bullet[\xi_1, \ldots \xi_n]$ can be turned into a differential bigraded algebra $(B, d_{k,\pm})$ with homology isomorphic to $\text{MSS}^\pm(\bigcirc_{cc})(\mathbb{C})$ where $\text{col} = V_\pm(n\omega_1)$. The grading on $B$ and the differential $d_{k,\pm}$ depends on $k$ and the sign $\pm$.*

**Remark 2.59.** A conjectural grading and differential is formulated in [60] for $-$. In case $k = n = 2$, Conjecture 2.58 follows up to an overall grading shift by a comparison of [34] with the formulas in [60], see [138] for a precise statement. In general the conjecture is open.

**Q5**: *Is there a conceptual method to compute the categorified coloured invariants?*
**Q6**: *To which extent is $\mathrm{MSS}^{\pm}_{\mathrm{col}}$ and its extension to framed tangles functorial?*

Motivated by and based on constructions of link homologies in physics, invariants of 3-manifolds are developed in e.g. **[62]**, **[63]**. Mathematically, steps in this direction are done in **[54]** by constructing categorified $3j$- and $6j$-symbols via fractional Euler characteristics.
**Q7**: *Do these coloured $\mathfrak{sl}_k$-invariants give rise to some invariant of 3-manifolds?*

### 3. Two proposals toward 4-TQFTs

We sketch two promising routes towards 4-TQFT based on Soergel bimodules.

**Braided monoidal structure on 2-representations.** Recall the starting point of algebraic categorification: the proposal **[35]** for constructing a 4-dimensional TQFT via Hopf categories. We like to interpret this as the wish of constructing, via categorified representation theory of quantum groups, a 0-1-2-3-4-*theory* **[52]**, i.e. a theory for $d = 4$ which not only evaluates at $d$- and $(d-1)-$, but also at $(d-2)$-, ..., 1- and 0-dimensional manifolds. To express the gluing laws between these levels one has to work **[8]**, **[94]** in general with an $n$-category of bordisms (viewed as $(\infty, n)$-category) and define a *fully extended n-TQFT* as a functor from this symmetric monoidal category into some symmetric monoidal $n$- (respectively $(\infty, n)$-)category. According to the *cobordism hypothesis* **[8]**, **[94]** such a fully extended TQFT $F$ is determined by the value $F(\mathrm{pt})$ at a point, see **[129]**, **[94]**, **[7]** for partial proofs.

Already the question *what does Chern–Simons theory attach to a point* is subtle and depends on the perspective. Following **[52]**, **[143]**, Chern–Simons theory or the related Witten–Reshetikhin–Turaev theory can be viewed as an *anomalous* 0-1-2-3 *theory* of oriented 4-manifolds, i.e. a morphism from the trivial theory to an invertible fully extended 4-TQFT $F$ defined on oriented manifolds. A similar interpretation was proposed by Walker, and the related invariant of a 4-manifold was combinatorially described in **[36]**. These interpretations propose to attach a certain *braided monoidal category* $F(\mathrm{pt})$ to a point **[52]**. Coming back to our setting, this suggests that a *categorification of the braided monoidal category of representations of a quantum group* might arise as the value $F_{\mathrm{cat}}(\mathrm{pt})$ of a point of an anomaly $F_{\mathrm{cat}}$, some fully extended (possibly partial) 5-TQFT with an anomalous 0-1-2-3-4-theory.

**Remark 3.1.** Some relevance **[102]** for 4-dimensional topology is already visible in $\mathrm{Kh}$ and $\mathrm{MSS}^-$, i.e. in categorified intertwiners of $\mathcal{F}und_2^-$ as in Remark 2.32, in particular via tangle cobordisms for surfaces in dimension 4 **[78]**, **[68]** and for invariants of 4-manifolds **[107]**.

Concretely, one seeks a monoidal structure on the 2-category of 2-representations of $U_v(\mathfrak{q})$ as in Remark 2.47 for $\mathfrak{q} = \mathfrak{gl}_k$ or $\mathfrak{sl}_k$ and say $\eta = 1$. Sections 2.3 and 2.4 presented tensor product categorifications and indicated categorifications of the duals and the braiding morphisms. The *process of taking tensor products*, i.e. the construction of a *tensor product* or an inner hom for 2-representations is however more involved. Inspired by (bordered)

Heegard-Floer theory, Manion and Rouquier [100], [101] give such a construction in case q is the positive part $\mathfrak{gl}(1|1)_+$ of the Lie superalgebra $\mathfrak{gl}(1|1)$ (for the analogue of Remark 2.47 see [76]). The passage to $\mathfrak{gl}(1|1)_+$ surprisingly simplifies the situation. In contrast to $\mathfrak{q} = \mathfrak{gl}_2$ or $\mathfrak{sl}_2$, homotopical complications disappear. The result of [101] is supposed to connect (as the value at an interval) to a slightly different type of TQFT and the theory predicted in [63].

**Remark 3.2.** This seemingly very different $\mathfrak{gl}(1|1)_+$-theory is still related to Section 2.3 via an interpretation in terms of subquotients of category $O$ [88]. Only $\mathfrak{gl}(1|1)_+$ appears, since categorical actions of $\mathfrak{gl}(1|1)$ have not yet been defined. This might be connected with the non-semisimplicity of the finite-dimensional representation theory of $\mathfrak{gl}(1|1)$, see e.g. [25].

Rouquier however announced (in an appropriate $A_\infty$-setting) the existence of a monoidal structure on the 2-category of 2-representations of $U_v(\mathfrak{q})$ for an arbitrary Kac–Moody Lie algebra q and a candidate for a braiding. This result should provide the desired value $F_{\text{cat}}(\text{pt})$. In the spirit of [35] we *propose to call* the resulting 2-category with its braided monoidal structure *the Hopf category of* q and reformulate ideas from [35] as:

**Prediction 3.3.** *The Hopf category of* q *is the value* $F_{\text{cat}}(\text{pt})$ *for an anomaly fully extended (partially defined at the top) 5-TQFT with an anomalous 0-1-2-3-4-theory.*

**Soergel bimodules, braided monoidal 2-categories and TQFT.** We finish by proposing another approach towards 4-TQFTs using more directly categories of Soergel bimodules. This is again motivated by the idea that a braided monoidal category might occur as the value $F(\text{pt})$ at a point [52] in a 0-1-2-3-theory. We seek to increase the dimensions to a 0-1-2-3-4-theory with a *braided monoidal bicategory* as the value $F(\text{pt})$ of some fully extended 5-TQFT $F$. We sketch some first steps. This is current work with Paul Wedrich.

**Remark 3.4.** The first definition of a *semistrict monoidal* and a *semistrict braided monoidal* 2-*category* is due to [71], [72]. It was then improved and put into a more concise definition in [9] with a technical adjustment in [37]. The concepts also appear as (braided) Gray monoids in [38]. By a *braided monoidal bicategory* we mean the less strict version from [64].

In the following let $m, n \in \mathbb{N}_0$. Recall the category $\mathcal{SB}im_n$ of Soergel bimodules from Section 2.2 with $R_0 := \mathbb{C}$ and $\mathcal{SB}im_0$ finite dimensional graded vector spaces. We view $\mathcal{SB}im_n$ as a graded monoidal category with tensor product $\circ_1 : (M, N) \mapsto M \otimes_{R_n} N$. If now $M \in \mathcal{SB}im_m$, $N \in \mathcal{SB}im_n$, then $M \boxtimes N := M \otimes_\mathbb{C} N$ is an $R_m \otimes_\mathbb{C} R_n = R_{m+n}$-bimodule and by construction an object in $\mathcal{SB}im_{m+n}$. For morphisms $f$ and $g$ in $\mathcal{SB}im_m$ and $\mathcal{SB}im_n$ respectively we define then $f \boxtimes g$ in the obvious way and set $m \boxtimes n = m + n$.

To get the desired semistrictness we use the monoidal category $\mathcal{DB}im_n \simeq \mathcal{SB}im_n$ [46], [48] from Remark 2.20. We omit giving the definition of $\mathcal{DB}im_n$ (it would not even fit on a page) and just recall that $\mathcal{DB}im_n$ is the Karoubian closure of a graded monoidal category $\mathcal{DB}im'_n$ [48]. The definition of $\mathcal{DB}im'_n$ is via generators and relations in terms of diagrams (similar to the usual string diagrams for higher categories). The morphism spaces come with distinguished bases, often called light leaves bases. A picked basis allows one to mimic the concept of coordinatised vector spaces from [71] and (semi)strictify the setup. Implicitly we assume this now, not altering the notation. We obtain categories *enriched in* $\mathbb{C}$*-linear categories* (we use *bicategories* as in [14] and *monoidal bicategories* as e.g. in [129]):

**Theorem 3.5.** *There is a bicategory* $^{(2)}\mathcal{SB}im$ *with objects* $\mathbb{N}_0$ *and nontrivial* $\mathbb{C}$*-linear morphism categories* $^{(2)}\mathcal{SB}im(m,n)$ *only in case* $m = n$, *in which case* $^{(2)}\mathcal{SB}im(n,n) = \mathcal{SB}im_n$ *with composition* $\circ_1$. *Similarly for* $^{(2)}\mathcal{DB}im$, *but with* $\mathcal{SB}im_n$ *replaced by* $\mathcal{DB}im_n$. *Moreover,* $^{(2)}\mathcal{SB}im$ *and* $^{(2)}\mathcal{DB}im$ *can be turned into monoidal bicategories with tensor functor* $\boxtimes$, *even into a semistrict monoidal* 2-*category in the sense of* [9], [37] *in case of* $^{(2)}\mathcal{DB}im$.

**Remark 3.6.** An analogue of $^{(2)}\mathcal{SB}im$ for singular Soergel bimodules as in Remark 2.23 exists as well (with the expected definition). For simplicity, we do not discuss this here.

The proof is done by *explicitly* constructing the required data and checking the coherence relations. Replacing $\mathcal{SB}im_n$ with the graded dg-category $C^b(\mathcal{SB}im_n)$ of bounded chain complexes of Soergel bimodules we get a category *enriched in graded* $\mathbb{C}$*-linear dg-categories*, similarly with $C^b(\mathcal{DB}im_n)$ instead of $\mathcal{DB}im_n$. Theorem 3.5 directly extends and provides bicategories, denoted $^2\mathcal{SB}im$ and $^2\mathcal{DB}im$, now realised as categories *enriched* [56] *in the monoidal category of* $\mathbb{C}$*-linear dg-categories* [75].

We consider from now on only the stricter version $^2\mathcal{DB}im$. To define a braiding, we need in particular an adjoint equivalence $\mathbb{B} : \boxtimes \Rightarrow \boxtimes^{op}$, [9], [64]. This data includes a *braiding* 1-*morphism* $\mathbb{B}((a,b))$ in $^2\mathcal{DB}im(a \boxtimes b = a + b, b \boxtimes a = b + a)$ for any $a, b \in \mathbb{N}_0$. Thinking intuitively about this braiding 1-morphism gives us a candidate:

(3.1) $\quad \beta = \underbrace{\phantom{XXX}}_{a}\underbrace{\phantom{XXX}}_{b} \rightsquigarrow \quad \ddot{\beta} = (\beta_b \cdots \beta_1) \cdots (\beta_{i+b-1} \cdots \beta_i) \cdots (\beta_{a+b-1} \cdots \beta_a)$

$\quad$ Rouquier complex $X(\beta)$ with

namely the Rouquier complex $X(\ddot{\beta}) \in C^b(\mathcal{SB}im_{a+b})$ from Section 2.2 with $\ddot{\beta}$ as in (3.1) translated via the above equivalence to an object $\mathbb{B}((a,b))$ in $C^b(\mathcal{DB}im_{a+b})$.

**Theorem 3.7.** *The proposed adjoint equivalence* $\mathbb{B}$ *satisfies the required naturality conditions* [91] *for the generating* 1- *and* 2-*morphisms of* $\mathcal{DB}im$ *up to canonical homotopy*.

To obtain however an honest braiding one has to pass to the homotopy categories which loses quite a lot of information or to a category $^2_\infty\mathcal{DB}im$ enriched in $\infty$-*categories* [56]. We construct such a category $^2_\infty\mathcal{DB}im$ by applying a (rather technical and not standard) dg-nerve construction to the morphism categories. We expect this construction to satisfy:

**Conjecture 3.8.** $^2_\infty\mathcal{DB}im$ *is a braided monoidal bicategory*.

**Remark 3.9.** Braided monoidal 2-categories with linear hom-categories and finiteness conditions should be objects in some symmetric monoidal 5-category which arises as next step in the ladder of symmetric monoidal $n$-categories (made explicit in [20]: objects are certain monoidal categories for $n = 3$ and certain braided monoidal categories for $n = 4$).

**Remark 3.10.** We can view $_\infty^2\mathcal{DB}im$ as a category object in $\infty$-categories $\mathbf{Cat}_\infty$. We expect this to be an $E_2$-algebra in the $\infty$-category of $(\infty,2)$-categories [65], [93]. Higher Morita theory of $E_n$-algebras [69] provides a possible ambient $(\infty,5)$-category for our hoped for TQFT.

Because of lacking finiteness conditions one should not expect $n$-dualisability [94] of $_\infty^2\mathcal{DB}im$ for $n > 3$, but we hope it holds for $n = 3, 4$ for quotients arising from actions on the $\mathfrak{gl}_k$-theories MSS_ for fixed $k \in \mathbb{N}$: An analogue of $_\infty^2\mathcal{DB}im$ defined using singular Soergel bimodules, Remark 3.6, acts by Remark 2.24 on the 2-categories $\hat{\mathcal{O}}^\mathbf{d}$ from Remark 2.47 for any fixed $k$. We conjecture that the largest quotient $_\infty^2\mathcal{DB}im(k)$ which still acts (for fixed $k$) has the desired finiteness properties to provide a fully extended (partial) 5-TQFT:

**Conjecture 3.11.** Soergel bimodules give rise to a braided monoidal bicategory $_\infty^2\mathcal{DB}im(k)$, $k \in \mathbb{N}$, which is the value at a point of an anomaly with an anomalous 0-1-2-3-4-theory.

### Acknowledgements
It is a pleasure to thank J. Brundan, A. Mathas, R. Rouquier, M. Stroppel and J. Sussan for many mathematical and (non-)mathematical discussions and for openly sharing ideas over the years. I am grateful to the Bonn representation theory group, in particular J. Eberhardt, G. Jasso, T. Heidersdorf, J. Matherne, J. Meinel, D. Tubbenhauer, P. Wedrich and T. Wehrhan for feedback and constructive criticism on draft versions of this article.

### Funding
This work was supported by the Hausdorff Center of Mathematics (HCM) in Bonn.

**Catharina Stroppel**

Mathematical Institute, Endenicher Allee 60, 53115 Bonn, stroppel@math.uni-bonn.de